\documentclass[11pt]{article}
\usepackage{amsmath,amsfonts,amssymb,amsthm}
\usepackage{xcolor}
\usepackage{color}
\usepackage{graphicx}
\AtBeginDocument{
	\label{CorrectFirstPageLabel}
	
}

\makeatletter
\newcommand\xleftrightarrow[2][]{%
	\ext@arrow 9999{\longleftrightarrowfill@}{#1}{#2}}
\newcommand\longleftrightarrowfill@{%
	\arrowfill@\leftarrow\relbar\rightarrow}
\makeatother

\usepackage{comment}

\theoremstyle{plain}

\newtheorem{theorem}{Theorem}[section]
\newtheorem{lemma}[theorem]{Lemma}
\theoremstyle{remark}
\newtheorem{definition}[theorem]{Definition}
\newtheorem{example}[theorem]{Example}

\usepackage{dsfont}
\def\E{{\mathbb E}}
\def\R{{\mathbb R}}
\def\P{{\mathbb P}}
\def\N{{\mathbb N}}

\def\L{{\lambda}}
\def\F{{\mathcal F}}

\def\vt{{\vartheta}}
\def\eps{{\varepsilon}}
\def\emptyset{{\varnothing}}

		\numberwithin{equation}{section}

\begin{document} 
\renewcommand{\thefootnote}{\fnsymbol{footnote}}
\author{{\sc Moritz Otto\footnotemark[1]}}
\footnotetext[1]{Department of Mathematics,
	Aarhus University, Aarhus, Denmark, otto@math.au.dk}

\date{}

\title{Poisson approximation of Poisson-driven point\\ processes and extreme values in stochastic geometry}

\maketitle

\begin{abstract}
	\noindent  We study point processes that consist of certain centers of point tuples of an underlying Poisson process. Such processes arise in stochastic geometry in the study of exceedances of various functionals describing geometric properties of the Poisson process. We use a coupling of the point process with its Palm version to prove a general Poisson limit theorem. We then combine our general result with the theory of asymptotic shapes of large cells (Kendall's problem) in random mosaics and prove Poisson limit theorems for large cells (with respect to a general size functional) in the Poisson-Voronoi and -Delaunay mosaic. As a consequence, we establish Gumbel limits for the asymptotic distribution of concrete size functionals and specify the rate of convergence. This extends extreme value results from Calka and Chenavier (2014) and Chenavier (2014).
\end{abstract}


\section{Introduction}\label{s1}
Point processes and random mosaics are fundamental objects in modern probability and find many applications, both in theory and practice. Often one is interested in certain geometric features of the configuration (e.g. the distance of a point to its nearest neighbor or the volume of a cell in the mosaic). One way to study these properties is to thin the process to a subprocess that is given by all points with a certain property (e.g. a large distance to its nearest neighbor). More generally, one can define a point process of certain centers of point tuples. In this article we study those processes in the situation where the original process is Poisson.

Let $\eta$ be a Poisson process on the locally compact second countable Hausdorff space $(\mathbb{X},\mathcal{X})$. We consider $\eta$ as a random object in the space $\mathbf{N}$ of locally finite counting measures on $\mathbb X$ equipped with its  standard $\sigma$-field $\mathcal{N}$ (see Section \ref{s2} for its definition). For $m\in \{1,\dots,d+1\}$, let $g:\mathbb{X}^m\times\mathbf{N}\to \{0,1\}$ be measurable and symmetric in the first $m$ coordinates and let $z:\mathbb{X}^m \to \mathbb{X}$ be a measurable and symmetric function. For instance, if the intensity measure of $\eta$ is absolutely continuous (such that all $\mathbf{x}\in \eta^{(m)}$ are a.s. in general position), we can choose $z(\mathbf{x})$ as the center of the unique $(m-2)$-sphere through $x_1,\dots,x_m$. We denote the Dirac measure in a point $z\in \mathbb X$ by $\delta_z$ and consider the point process
\begin{align}
	\xi:=\frac{1}{m!} \sum\limits_{\mathbf{x}\in \eta^{(m)}}  g(\mathbf{x},\eta)\,\delta_{z(\mathbf{x})}. \label{inxidef}
\end{align}
Here, the function $g$ can be understood as a selection mechanism that decides if $z(\mathbf x)$ is considered or not. We will require that $g(\mathbf{x},\mu)$ depends only locally on the configuration $\mu$ around $z(\mathbf{x})$ and formalize this concept in Definition \ref{defstab}. In many applications, there is some compact $W \subset \R^d$ and an underlying function $f:\mathbb{X}^m \times \mathbf{N}\to \R$ that measures the size of a geometrically defined object associated to $\mathbf{x}$ (e.g. its Voronoi or Delaunay cell with respect to $\mu$) and $g(\mathbf{x},\mu):=\mathbf{1}\{z(\mathbf{x}) \in W\}\mathbf{1}\{f(\mathbf{x},\mu)>v\}$ is chosen as the indicator that encodes whether $z(\mathbf{x}) \in W$ and whether $f(\mathbf{x},\mu)$ exceeds a certain value $v \in \R$ or not.

The motivation of our paper is twofold. First, we establish a general Poisson process approximation result for the process $\xi$ under the asumption that $g$ satisfies a stabilization condition. Our method of proof uses Stein's method and constructs a concrete coupling of $\xi$ and its Palm measure. Thus, we avoid the use of Glauber dynamics on the Poisson space which is applied in \cite{bobrowski2021poisson} to establish Poisson approximation. Second, we demonstrate the flexibility of our approach and derive Poisson approximation results for large cells in the Poisson-Voronoi and -Delaunay mosaics for general size functionals. This extends and generalizes results from \cite{calka2014extreme} and \cite{chenavier2014general}. We make use the theory of the asymptotic shape of large typical cells in random mosaics (Kendall's problem) that is contained in a series articles (see  \cite{{kovalenko1997proof},{hug2004large},{hug2004limit},{hug2005large},{hug2004large1},{hug2007asymptotic},{hug2007typical},{bonnet2018cells}}). 

Various results on Poisson approximation of point processes in stochastic geometry can be found in the literature. In \cite{schuhmacher2005distance} and \cite{schuhmacher2009distance} approximation results for dependent thinnings of point procesesses that have a density with respect to a Poisson process are derived. The Malliavin calculus on the Poisson space is used in \cite{schulte2012scaling} and \cite{schulte2016poisson} to determine scaling limits for Poisson U-statistics. In \cite{decreusefond2016functional} the theory of Glauber dynamics for birth and death processes is combined with the Chen-Stein method and a Poisson approximation result in Kantorovich-Rubinstein distance is derived for U-statistics. \cite{chenavier2014general} uses a d Poisson approximation result from \cite{arratia1989two} and \cite{arratia1990poisson} together with a discretizetion technique to study extremes in the Poisson-Voronoi and -Delaunay mosaic. \cite{chenavier2016extremes} and \cite{otto2021extremal} derive Poisson approximation results for small and large cells in the Poisson hyperplane mosaic. Poisson approximation for the volume of $k$-nearest neighbor balls is discussed in \cite{gyorfi2018limit}, \cite{chenavier2022limit}, \cite{bobrowski2021poisson} in the Euclidean space and in \cite{otto2022large} in the hyperbolic space. Moreover, Poisson approximation for functionals of Gibbs process is discussed in \cite{{schuhmacher2005distance},{schuhmacher2009distance},{last2023+disagreement}}. 
For more information on extreme values in random tessellations we refer to \cite{{bonnet2020maximal},{chenavier2019largest},{chenavier2018cluster},{chenavier2015extremal},{bonnet2018small},{schneider2019small}}.

This article is organized as follows. We settle our notation and provide background information on Palm theory and stopping sets in Section 2. In Section 3 we introduce the notion of stabilization that we work with and state and prove a general theorem on Poisson approximation for point processes. In Section 4 we study cells in the Poisson-Voronoi mosaic that are large with respect to a general size functional. Beyond that, we consider large Poisson-Delaunay cells in Section 5. 
\section{Preliminaries}\label{s2}
\subsection{Point and Poisson processes}
Let $(\mathbb{X},\mathcal X)$ be a locally compact second countable Hausdorff space. Let $\mathbf{N}$ denote the space of locally finite (i.e. finite on compact subsets) counting measures $\mu$ on $\mathbb{X}$ and let $\mathcal{N}$ be the $\sigma$-field on $\mathbf{N}$ generated by the mappings $\mu \mapsto \mu(B)$, $B \in \mathcal X$. Examples of elements of $\mathbf{N}$ are the zero measure $0$ and the Dirac measure $\delta_x$ in the point $x \in \mathbb X$, given by $\delta_x(B):=\mathbf{1}_B(x),\,B \in \mathcal{X}$. We write $\mu_B:=\mu \cap B$ for the restriction of $\mu \in \mathbf N$ to $B \in \mathcal X$. For $k \in \N$, $\mathbf{x}=(x_1,\dots,x_k)\in \mathbb{X}^k$, $\mu \in \textbf{N}$ and $A \in \mathcal{X}$ we write $\delta_{\mathbf{x}}:=\delta_{x_1}+\cdots+\delta_{x_k}$, $\mu_{\mathbf{x}}:=\mu+\delta_{\mathbf{x}}$ and $\mathbf{x} \in A$ if $x_i\in A$ for every $i \in [k]$.  By a slight abuse of notation, we write $x \in \mu$ if the point $x$ is charged by the configuration $\mu$, i.e. $\mu(\{x\})>0$. 

Suppose $\mu \in \mathbf{N}$ is given by $\mu=\sum_{i=1}^k \delta_{x_i}$ for some $k \in \N_0 \cup \{\infty\}$ and some $x_1,\dots,x_k \in \mathbb{X}$ (not necessarily distinct). For $m \in \N$ we define the {\em factorial measure} (see \cite{last2017lectures}, (4.5)) $\mu^{(m)}$ of $\mu$ by
\begin{align*}
	\mu^{(m)}:=\sum^{\neq}\limits_{i_1,\dots,i_m \le k} \delta_{(x_{i_1},\dots,x_{i_m})},
\end{align*}
where the superscript $\neq$ indicates summation over $m$-tuples with pairwise different entries and where an empty sum is defined as zero. (For $k=\infty$ this involves only integer-valued indices.) A measure $\mu \in \textbf{N}$ is called {\em simple} if $\mu(\{x\})\in \{0,1\}$ for all $x \in \mathbb{R}^d$. Let $\mathbf{N}_s$ denote the set of simple locally finite counting measures. It is a measurable subset of $\mathbf{N}$ (see \cite{schneider2008stochastic}, p. 51) and its induced $\sigma$-field is denoted by $\mathcal{N}_s$.

A {\em point process} is a random element $\xi$ in $\textbf{N}$, defined over some fixed probability space $(\Omega, \mathcal{A},\P)$. We assume that this probability space is rich enough to support all random objects in this article. By definition, $\xi(B)$ is a random variable in $\N_0 \cup \{\infty\}$ for every $B \in \mathcal{X}$. A point process $\xi$ is called {\em simple} if $\P(\xi \text{ is simple})=1$. 

The central object of this article is the {\em Poisson process}. We refer to \cite[Chapter 3]{last2017lectures} for its definition and basic properties. Particularly useful is the following {\em multivariate Mecke equation} (see \cite[Theorem 4.5]{last2017lectures}):

Let $\eta$ be a Poisson process on $\mathbb{X}$ with $\sigma$-finite intensity measure $\lambda$. Then we have for all $k \in \N$ and every measurable function $f:\mathbb{X}^k\times \mathbf{N}\to [0,\infty]$,
\begin{align}
	\E \int f(\mathbf{x},\eta )\,\eta^{(k)}(\mathrm{d}\mathbf{x})=\int \E f(\mathbf{x},\eta+\delta_{\mathbf x})\,\lambda^k(\mathrm{d}\mathbf{x}). \label{mecke}
\end{align}

\subsection{Palm measures and stopping sets}
Following \cite{kallenberg2017random}, Chapter 6, we next introduce Palm processes. Let $\eta,\xi$ be two point processes and assume that $\xi$ has $\sigma$-finite intensity measure $\E \xi$. In general there exists a whole family of Palm measures, one for every $x \in \mathbb{R}^d$. Palm measures generalize the notion of regular conditional distributions (see Theorem 6.3 in \cite{kallenberg2006foundations}) and agree with them for $\xi=\delta_x$ for some $x \in \mathbb{R}^d$. The $\sigma$-finiteness of $\E \xi$ implies that the {\em Campbell measure} $C_{\eta,\xi}$ defined by
\begin{align*}
	C_{\eta,\xi}f:= \E \int f(x,\eta)\,\xi(\mathrm{d}x),\quad f:\R^d \times \mathbf{N} \to [0,\infty) \text{ measurable},
\end{align*}
is also $\sigma$-finite. There exists a (unique) probability kernel $P_{\eta,\xi}^x$ from $\R^d$ to $\textbf{N}$ such that for all measurable $f:\R^d \times \mathbf{N} \to [0,\infty)$ the disintegration
\begin{align*}
	C_{\eta,\xi}f= \iint f(x,\mu)\,P_{\eta,\xi}^x(\mathrm{d}\mu)\,\E \xi(\mathrm{d}x)
\end{align*}
holds. The measure $P_{\eta,\xi}^x$ is called the {\em Palm measure} of $\eta$ with respect to $\xi$ at $x$. A point process $\eta^{\xi,x}$ with distribution $P_{\eta,\xi}^x$ is called a {\em Palm process} or a {\em Palm version} of $\eta$ with respect to $\xi$ at $x$. If $\xi$ is simple, $\eta^{\xi,x}$  can be interpreted as the process $\eta$ seen from $x$ conditioned on $\xi$ having a point in $x$. If $\xi = \eta$ a.s.\ we write $\eta^x$ for a Palm version of $\eta$ (with respect to itself) at $x$ and we obtain from Lemma 6.2(ii) in \cite{kallenberg2017random} that
\begin{align*}
	\P(x \in \eta^x)=1.
\end{align*}
The process $\eta^x-\delta_x$ is called a {\em reduced Palm process} of $\eta$ at $x$.

An important tool in the analysis of point processes are stopping sets. They generalize the concept of stopping times for random variables.  Let $\mathcal{F}$ denote the system of closed sets in $\mathbb X$. We endow $\mathcal{F}$ with the smallest $\sigma$-field containing $\mathcal{F}_K:=\{F \in \mathcal{F}:\,F \cap K \neq \emptyset\}$ for all compact $K \subset \mathbb X$. For $F \in \mathcal{F}$ we denote by $\mu_F$ the restriction of $\mu$ to $F$. Moreover, let $\mathcal{N}_F$ be the $\sigma$-field on $\mathbf{N}$ generated by the mappings $\mu\mapsto \mu(B \cap F)$, $B \in \mathcal{X}$.
\begin{definition} 
	A measurable map $S: \textbf{N} \to \mathcal{F}$ is called {\em stopping set} (with respect to the filtration $(\mathcal{N}_F)_{F \in \mathcal{F}}$) if $\{\mu \in \textbf{N}:\,S(\mu)\subset F\} \in \mathcal{N}_F$ for all $F \in \mathcal{F}$.
\end{definition}

Intuitively, if $S$ is a stopping set and $\eta$ is a random element in $\mathbf{N}$, then $S(\eta)$ is a random subset of $\mathbb{R}^d$ such that $S(\eta)$ only depends on the restriction $\eta_{S(\eta)}$ of $\eta$ to $S(\eta)$.

From \cite[Proposition A.1]{baumstark2009gamma} we have that a measurable map $S:\mathbf{N}\to \F$ is a stopping set if and only if $S(\mu)=S(\mu_{S(\mu)})$ for all $\mu \in \mathbf{N}$ and if the following implication holds for all $\mu,\varphi \in \mathbf{N}$:
\begin{align}
	\varphi=\mu_{S(\varphi)} \Longrightarrow S(\varphi)=S(\mu). \label{propstop}
\end{align} 

The following lemma is similar to Lemma A.2 in \cite{last2021phase}. Since the latter works under a more general notion of stopping sets, we give a proof of the statement.
\begin{lemma} \label{remstop} Let $\varphi, \mu \in \mathbf{N}$ and let $S:\mathbf{N}\to \mathcal{F}$ be a stopping set. Then it holds that
	\begin{align*}
		\varphi=	\mu_{S(\mu)} \quad \Longleftrightarrow \quad \varphi= \mu_{S(\varphi)}
	\end{align*}
	and in this case we also have that $S(\varphi)=S(\mu)$.
\end{lemma}

\begin{proof}
	Let $S:\mathbf{N}\to \mathcal{F}$ be a stopping set and assume that $\varphi=\mu_{S(\mu)}$. Then we have by \cite[Proposition A.1]{baumstark2009gamma} that $S(\mu)=S(\mu_{S(\mu)})=S(\varphi)$. For the other direction we assume that $\varphi=\mu_{S(\varphi)}$. Then it follows from \eqref{propstop} that $S(\varphi)=S(\mu)$ and, hence, that $\varphi=\mu_{S(\mu)}$.
\end{proof}

The next statement will be used repeatedly in Section 3. 

\begin{lemma} \label{stopmecam}
	Let $\eta$ and $\eta'$ be independent Poisson process on $\mathbb{X}$ with $\sigma$-finite intensity measure $\lambda$, let $S$ be a stopping set such that $S(\eta)$ is a.s.\ compact and let $h:\mathbf{N}\times\mathbf{N}\to [0,\infty)$ be measurable. Then we have
	\begin{align*}
		\E h(\eta_{S(\eta)},\eta_{S(\eta)^c})=\E h(\eta_{S(\eta)},\eta'_{S(\eta)^c}).
	\end{align*}
\end{lemma}

\begin{proof}
	The argument can be assembled from different sources in the literature (see \cite{baumstark2009gamma}, \cite{zuyev1999stopping}). Nevertheless, we give a proof for completeness and convenience of the reader. We have
	\begin{align}
		\E h(\eta_{S(\eta)},\eta_{S(\eta)^c})&= \sum\limits_{k=0}^\infty \E h(\eta_{S(\eta)},\eta_{S(\eta)^c})\,\mathbf{1}\{\eta(S(\eta))=k\}\nonumber\\
		&= \sum\limits_{k=0}^\infty \frac{1}{k!} \E \int h(\delta_{\mathbf{x}},\eta_{S(\eta)^c})\,\mathbf{1}\{\eta_{S(\eta)}=\delta_{\mathbf{x}}\}\,\eta^{(k)}(\mathrm{d}\mathbf{x}),\label{funstdis}
	\end{align}
	where we recall that $\delta_{\mathbf{x}}:=\delta_{x_1}+\cdots+\delta_{x_k}$ for $\mathbf{x}=(x_1,\dots,x_k)\in \mathbb X^k$. By Lemma \ref{remstop} we have that  $\eta_{S(\eta)}=\delta_{\mathbf{x}}$ if and only if $\eta_{S(\delta_{\mathbf{x}})}=\delta_{\mathbf{x}}$ and in this case it holds that $S(\eta)= S(\delta_{\mathbf{x}})$. Hence, \eqref{funstdis} is given by
	\begin{align*}
		\sum\limits_{k=0}^\infty \frac{1}{k!}\E \int h(\delta_{\mathbf{x}},\eta_{S(\delta_{\mathbf{x}})^c})\,\mathbf{1}\{\eta_{S(\delta_{\mathbf{x}})}=\delta_{\mathbf{x}}\}\,\eta^{(k)}(\mathrm{d}\mathbf{x}).
	\end{align*}
	From the multivariate Mecke equation \eqref{mecke} we obtain that the above is given by
	\begin{align}
		\sum\limits_{k=0}^\infty \frac{1}{k!} \int\E h(\delta_{\mathbf{x}},(\eta+\delta_{\mathbf{x}})_{S(\delta_{\mathbf{x}})^c})\,\mathbf{1}\{(\eta+\delta_{\mathbf{x}})_{S(\delta_{\mathbf{x}})}=\delta_{\mathbf{x}}\}\,\lambda^{k}(\mathrm{d}\mathbf{x}).\label{funstdism}
	\end{align}
	Note that $(\eta+\delta_{\mathbf{x}})_{S(\delta_{\mathbf{x}})}=\delta_{\mathbf{x}}$ implies that $(\eta+\delta_{\mathbf{x}})_{S(\delta_{\mathbf{x}})^c}=\eta_{S(\delta_{\mathbf{x}})^c}$. Let $\eta'$ be a point process that is independent of $\eta$ with $\eta \stackrel{d}{=}\eta'$. Since $\eta_{S(\delta_{\mathbf{x}})}$ and $\eta_{S(\delta_{\mathbf{x}})^c}$ are independent, \eqref{funstdism} is given by
	\begin{align*}
		&\sum\limits_{k=0}^\infty \frac{1}{k!} \int\E h(\delta_{\mathbf{x}},\eta'_{S(\delta_{\mathbf{x}})^c})\,\mathbf{1}\{(\eta+\delta_{\mathbf{x}})_{S(\delta_{\mathbf{x}})}=\delta_{\mathbf{x}}\}\,\lambda^{k}(\mathrm{d}\mathbf{x})\\
		&\quad =\sum\limits_{k=0}^\infty \frac{1}{k!} \E\int h(\delta_{\mathbf{x}},\eta'_{S(\delta_{\mathbf{x}})^c})\,\mathbf{1}\{\eta_{S(\delta_{\mathbf{x}})}=\delta_{\mathbf{x}}\}\,\eta^{(k)}(\mathrm{d}\mathbf{x}),
	\end{align*}
	where we have applied the Mecke equation to obtain the equality. Using here Lemma \ref{remstop} again, we arrive at 
	\begin{align*}
		&\sum\limits_{k=0}^\infty \frac{1}{k!} \E \int h(\delta_{\mathbf{x}},\eta'_{S(\eta)^c})\,\mathbf{1}\{\eta_{S(\eta)}=\delta_{\mathbf{x}}\}\,\eta^{(k)}(\mathrm{d}\mathbf{x})\nonumber\\
		&\quad=	\sum\limits_{k=0}^\infty \E h(\eta_{S(\eta)},\eta'_{S(\eta)^c})\,\mathbf{1}\{\eta(S(\eta))=k\}=\E h(\eta_{S(\eta)},\eta'_{S(\eta)^c}).
	\end{align*}
\end{proof}

\section{General result on Poisson process approximation}\label{s3}

Let $m \in [d+1]$ and let $\eta$ be a Poisson process in $\mathbb X$ with $\sigma$-finite and diffuse intensity measure $\lambda$. In the following we assume that $g:\mathbb X^m\times \mathbf{N}\to \{0,1\}$ is measurable and symmetric in the first $m$ coordinates, i.e. 
\begin{align*}
	g(x_1,\dots,x_m,\mu)=g(x_{\pi(1)},\dots,x_{\pi(m)},\mu)
\end{align*}
for all $x_1,\dots,x_m \in \mathbb X,\,\mu \in \mathbf{N}$ and every permutation $\pi:[m]\to [m]$. We think of $g$ as a selection mechanism that decides whether an $m$-tuple $\mathbf{x}=(x_1,\dots,x_m)$ is considered or not. Let $z:\mathbb{X}^m \to \mathbb{X}$ be a measurable and symmetric function.
We define the point process
\begin{align}
	\xi [\mu]=\frac{1}{m!} \sum\limits_{\mathbf{x}\in \mu^{(m)}}  g(\mathbf{x},\mu)\,\delta_{z(\mathbf{x})},\quad \mu \in \mathbf N, \label{xidef}
\end{align}
and write $\xi:=\xi[\eta]$. Note that for $m=1$ and $z(x)=x$, the process $\xi$ is a thinning of the Poisson process $\eta$. From the multivriate Mecke equation \ref{mecke} we find that the intensity measure $\E \xi$ of $\xi$ is given by
\begin{align}
	\E \xi(A)&= \frac {1}{m!} \E \int \mathbf{1}\{z(\mathbf{x}) \in A\}\,g(\mathbf{x},\eta)\,\eta^{(m)}(\mathrm{d}\mathbf{x})\nonumber\\
	&= \frac {1}{m!}\int \mathbf{1}\{z(\mathbf{x}) \in A\}\, \E g(\mathbf{x},\eta+\delta_\mathbf{x})\,\lambda^m(\mathrm{d}\mathbf{x}),\quad A \in \mathcal{X}.\label{intpalm}
\end{align}

The goal of this section is to approximate $\xi[\eta]$ by a Poisson process under the condition that $g$ is stabilizing. This concept is defined formally in the next definition. Loosely sopken, it requires that the value of $g(\mathbf{x},\mu)$ is determined by the resctriction of $\mu$ to a ball centred at $z(\mathbf x)$ with a finite radius.  We write $B_r(z)$ for the closed ball with radius $r>0$ around $z\in \mathbb{X}$.

\begin{definition} \label{defstab}
	Let $g:\mathbb{X}^m\times \mathbf{N}\to [0,\infty)$ be measurable and symmetric in the first $m$ coordinates and let $\eta$ be a Poisson process in $\mathbb{X}$ with $\sigma$-finite intensity measure. We call $g$ {\em stabilizing} if there exists a measurable function $R:\mathbb X\times \mathbf{N}\to [0,\infty]$, such that for all $\mathbf{x} \in \mathbb{X}^m$ we have
	\begin{enumerate}
		\item[(i)] $R(z(\mathbf{x}),\eta+\delta_\mathbf{x}) <\infty\quad \P$-a.s.
		\item[(ii)] $g(\mathbf{x},\mu)=g(\mathbf{x},\mu \cap B_{R(z(\mathbf{x}),\mu+\delta_\mathbf{x})}(z(\mathbf{x}))),\quad \mu \in \mathbf{N}.$
		\item[(iii)] The map $\mu \mapsto B_{R(z(\mathbf{x}),\mu+\delta_\mathbf{x})}(z(\mathbf{x}))$ from $\mathbf{N}$ to $\mathcal{F}$ is a stopping set.
		\item[(iv)] $\mathbf{x} \in B_{R(z(\mathbf{x}),\mu+\delta_\mathbf{x})}(z(\mathbf{x}))$,\quad $\mu\in \textbf{N}$.
	\end{enumerate}
	We call $R:= R(z,\mu)$ {\em stabilization radius} and use the notation $S(z,\mu):=B_{R(z,\mu)}(z),\,z \in \R^d,\,\mu \in \mathbf{N}$.  
\end{definition}

Note that from Definition \ref{defstab}(ii) and (iii) it follows that
\begin{align}
	g(\mathbf{x},\mu)=g(\mathbf{x},\mu \cap B_{R(z(\mathbf{x}),\mu+\delta_\mathbf{x})}(z(\mathbf{x}))+\chi \cap B_{R(z(\mathbf{x}),\mu+\delta_\mathbf{x})}(z(\mathbf{x}))^c),\quad \mu,\,\chi \in \mathbf N.\label{stoprem}
\end{align}
This property is sometimes assumed in the literature (see e.g. \cite{penrose2005normal}).

In order to state the main result of this section, we still need to fix some notation. For point processes $\xi$ and $\nu$ on $\mathbb{X}$ the {\em total variation distance} is given by 
\begin{align*}
	\mathbf{d_{TV}}(\xi,\nu):=\sup_{A \in \mathcal N} |\mathbb{P}(\xi \in A)-\mathbb{P}(\nu \in A)|.
\end{align*}
Moreover, we denote by $\mu_+$ and $\mu_-$ the positive and negative part of a finite signed measure $\mu$ on $\mathcal X$ and by $\|\mu\|:=\mu_+(\mathbb{X})+\mu_-(\mathbb{X})$ its total variation. 

%

\begin{theorem} \label{maithmallg}
	Let $\xi=\xi[\eta]$ be the process defined at \eqref{xidef} with stabilizing $g$ in the sense of Definition \ref{defstab} and assume that $\mathbb E \xi(\mathbb{X})<\infty$. Let $\nu$ be a Poisson process on $\mathbb{X}$ satisfying $\E \nu(\mathbb{X})<\infty$ and let $\mathbf x \mapsto b_{\mathbf x}$ be a measurable function from $\mathbb X^m$ to $[0,\infty)$. Then we have
	\begin{align}
		\mathbf{d_{TV}} (\xi, \,\nu)\le \|\E\xi-\E\nu\|+ T_{1}+T_2+T_3+T_4+T_5 
	\end{align}
	where 
	\begin{align*}
		&T_{1}:=\frac{2\mathbb{E}\xi(\mathbb{X})}{m!} \int_{\mathbb{X}^m} \mathbb{E}g(\mathbf{x},\eta+\delta_{\mathbf{x}})\mathbf 1\{R(z(\mathbf{x}),\eta+\delta_{\mathbf x})>b_{\mathbf{x}}\}  \,\lambda^m(\mathrm{d} \mathbf x),\\
		&T_{2}:=\frac{1}{(m!)^2}\int_{\mathbb{X}^m} \int_{\mathbb{X}^m} \mathbb{E} g(\mathbf{x},\eta+\delta_{\mathbf{x}})\mathbb{E}g(\mathbf{y},\eta+\delta_{\mathbf{y}})\,\mathbf 1 \{\|z(\mathbf{x})-z(\mathbf{y})\|\le b_{\mathbf{x}}+b_{\mathbf{y}}\}\,\lambda^m(\mathrm{d} \mathbf y)\,\lambda^m(\mathrm{d} \mathbf x),\\
		&T_{3}:=\frac{2}{(m!)^2}\int_{\mathbb{X}^m} \int_{\mathbb{X}^m}  \E g(\mathbf{x},\eta+\delta_{(\mathbf{x},\mathbf{y})})\,g(\mathbf{y},\eta+\delta_{(\mathbf{x},\mathbf{y})}) \,\mathbf{1}\{R(z(\mathbf{x}),\eta+\delta_{\mathbf{x},\mathbf{y}})>b_{\mathbf{x}}\}\, \lambda^m(\mathrm{d}\mathbf{y})\,\lambda^m(\mathrm{d}\mathbf{x}),\\
		&T_4:=\frac{1}{(m!)^2}\int_{\mathbb{X}^m} \int_{\mathbb{X}^m}  \E g(\mathbf{x},\eta+\delta_{(\mathbf{x},\mathbf{y})})\,g(\mathbf{y},\eta+\delta_{(\mathbf{x},\mathbf{y})}) \,\mathbf{1}\{\|z(\mathbf{x})-z(\mathbf{y})\|\le b_{\mathbf{x}}+b_{\mathbf{y}}\}\\
		&\qquad \qquad \times \mathbf{1}\{R(z(\mathbf{x}),\eta+\delta_{(\mathbf{x},\mathbf{y})})\le b_{\mathbf{x}}\}\,\mathbf{1}\{R(z(\mathbf{y}),\eta+\delta_{(\mathbf{x},\mathbf{y})})\le b_{\mathbf{y}}\}\, \,\lambda^m(\mathrm{d}\mathbf{y})\,\lambda^m(\mathrm{d}\mathbf{x}),\\
		&T_5:=\sum\limits_{\ell=1}^{m-1}\frac{1}{m!\,(m-\ell)!}\int_{\mathbb{X}^m}	\int_{\mathbb{X}^{m-\ell}} \mathbb{E} g((\mathbf{x}_\ell,\mathbf{y}),\eta+\delta_{(\mathbf{x},\mathbf{y})}) g(\mathbf x,\eta+\delta_{(\mathbf{x},\mathbf{y})})\lambda^{m-\ell}(\mathrm{d}\mathbf{x}) \lambda^{m} (\mathrm{d}\mathbf y)
	\end{align*}
	with $\mathbf{x}_\ell:=(x_1,\dots,x_\ell)$.
\end{theorem}

The strategy of the proof of Theorem \ref{maithmallg} is as follows. First, we construct a special (reduced) Palm version $\xi^{w!}$ of $\xi=\xi[\eta]$ at a given $w \in \mathbb{X}$. Thereafter, we find bounds on the total variations of the positive and the negative part of $\xi-\xi^{w!}$. Finally, we combine the two bounds and conclude the proof of Theorem \ref{maithmallg} using a general Poisson approximation result from \cite{barbour1992stein}. 

\begin{lemma}\label{lemxiz} Let $\xi=\xi[\eta]$ be the process defined at \eqref{xidef} and $\eta^{\xi,w}$ be a Palm process of $\eta$ with respect to $\xi$ at $w$ such that $\eta^{\xi,w}$ and $\eta$ are independent point processes. For $\mathbb E \xi$-almost all $w \in \mathbb{X}$,
	\begin{align}
		\xi^{w!}:=\xi[(\eta^{\xi,w})_{S(w,\eta^{\xi,w})}+\eta_{S(w,\eta^{\xi,w})^c}]-\delta_w \label{def:xiz}
	\end{align}
	is a reduced Palm process of $\xi$ (with respect to itself) at $w$.
\end{lemma}

\begin{proof}
	It suffices to show that
	\begin{align}
		(\eta^{\xi,w})_{S(w,\eta^{\xi,w})}+\eta_{S(w,\eta^{\xi,w})^c} \label{etazetaclaimallg}
	\end{align}
	is a Palm version of $\eta$ with respect to $\xi$ at $w$. To this end, let $h:\mathbb X\times \textbf{N} \to [0,+\infty)$ be measurable and let $\widetilde{\eta}$ be a Poisson process with $\eta \stackrel{d}{=} \widetilde{\eta}$ such that $\eta$ and $\widetilde{\eta}$ are independent. By independence of $\eta$ and $\eta^{\xi,w}$ and the definitions of $\eta^{\xi,w}$ and $\xi$, we obtain that
	\begin{align*}
		& \int \E h(w,(\eta^{\xi,w})_{S(w,\eta^{\xi,w})}+\eta_{S(w,\eta^{\xi,w})^c})\,\E \xi(\mathrm{d}w)\nonumber\\
		&\quad =\E \int h(w,\eta_{S(w,\eta)}+\widetilde{\eta}_{S(w,\eta)^c})\,\xi[\eta](\mathrm{d}w)\nonumber\\ 
		&\quad =\frac{1}{m!} \E \int  h(z(\mathbf{x}),\eta_{S(z(\mathbf{x}),\eta)}+\widetilde{\eta}_{S(z(\mathbf{x}),\eta)^c})\,g(\mathbf{x},\eta)\,\eta^{(m)} (\mathrm{d}\mathbf{x}),
	\end{align*}
	which is by the multivariate Mecke equation \eqref{mecke} given by
	\begin{align}
		\frac{1}{m!} \int	&\E  h(z(\mathbf{x}),\eta_{S(z(\mathbf{x}),\eta+\delta_\mathbf{x})}+\widetilde{\eta}_{S(z(\mathbf{x}),\eta+\delta_\mathbf{x})^c}+\delta_\mathbf{x})\,g(\mathbf{x},\eta+\delta_{\mathbf{x}})\,\lambda^m(\mathrm{d}\mathbf{x}),\label{etazetaallgmecke}
	\end{align}
	where we have used that $(\eta+\delta_{\mathbf x})_{S(z(\mathbf{x}),\eta+\delta_\mathbf{x})}=\eta_{S(z(\mathbf{x}),\eta+\delta_\mathbf{x})}+\delta_{\mathbf x}$ by Definition \ref{defstab}(iv). Now we apply Theorem \ref{stopmecam} with 	
	\begin{align*}
		g:\textbf{N}\times \textbf{N}\to [0,+\infty),\quad (\mu,\phi) \mapsto h(z(\mathbf{x}),\mu+\phi+\delta_\mathbf{x}) g(\mathbf{x},\mu+\delta_{\mathbf{x}}).
	\end{align*}
	Since $g$ is stabilizing and since $S$ is a stopping set, \eqref{etazetaallgmecke} can be written as
	\begin{align}	
		&\frac{1}{m!} \int \E h(z(\mathbf{x}),\eta+\delta_\mathbf{x})\,g(\mathbf{x},\eta+\delta_{\mathbf{x}})  \,\lambda^m(\mathrm{d}\mathbf{x}),
	\end{align}
	which is by the multivariate Mecke equation \eqref{mecke} and the definition of $\xi$ given by
	\begin{align*}
		&\frac{1}{m!} \E \int h(z(\mathbf{x}),\eta)\,g(\mathbf{x},\eta) \,\eta^{(m)}(\mathrm{d}\mathbf{x}) =\E \int h(w,\eta)\,\xi(\mathrm{d}w).
	\end{align*}
	Hence, $(\eta^{\xi,w})_{S(w,\eta^{\xi,w})}+\eta_{S(w,\eta^{\xi,w})^c}$ is indeed a Palm version of $\eta$ with respect to $\xi$ at $w$ for $\E \xi$-almost all $w$. 
\end{proof}

\begin{lemma}\label{lem+}
	Let $\xi=\xi[\eta]$ be the process defined at \eqref{xidef} and for $\mathbb{E} \xi$-almost all $w \in \mathbb{X}$ let $\xi^{w!}$ be the process defined at \eqref{def:xiz}. We have
	\begin{align*}
		&\int \mathbb{E} [(\xi-\xi^{w!})_+(\mathbb{X})]\,\E\xi(\mathrm{d}w)\le T_1+T_2,
	\end{align*}
	where $T_1, T_2$ are given in Theorem \ref{maithmallg}.
\end{lemma}

\begin{proof}
	Let $\eta$ and $\tilde \eta$ be independent Poisson processes on $\mathbb{X}$ with intensity measure $\lambda$. It follows from the definition of $\xi^{w!}$ in \eqref{def:xiz} that the left-hand side of the inequality in the lemma is given by
	\begin{align}
		\mathbb{E} \int(\xi[\eta]-\xi[\tilde \eta_{S(w,\tilde \eta)}+\eta _{S(w,\tilde \eta)^c}]+\delta_w)_+(\mathbb{X})\,\xi[\tilde \eta](\mathrm{d}w).\label{eqn:tv+lhs}
	\end{align}
	By definition of $\xi$, the integrand can be bounded as follows
	\begin{align}
		&(\xi[\eta]-\xi[\tilde \eta_{S(w,\tilde \eta)}+\eta _{S(w,\tilde \eta)^c}]+\delta_w)_+(\mathbb{X})\nonumber\\
		&\quad\le \frac {1}{m!}\int_{\mathbb{X}^m} g(\mathbf x,\eta) \mathbf 1\{\mathbf x \in (\tilde \eta_{S(w,\tilde \eta)}+\eta _{S(w,\tilde \eta)^c})^{(m)} \} \mathbf 1 \{g(\mathbf x,\eta)\neq g(\mathbf x,\tilde \eta_{S(w,\tilde \eta)}+\eta _{S(w,\tilde \eta)^c})\}\nonumber\\
		&\qquad \qquad \times \mathbf 1\{z(\mathbf x) \neq w\}\,\eta^{(m)}(\mathrm{d}\mathbf x)\label{eqn:gineq1}\\
		&\qquad +\frac {1}{m!}\int_{\mathbb{X}^m} g(\mathbf x,\eta) \mathbf 1\{\mathbf x \notin (\tilde \eta_{S(w,\tilde \eta)}+\eta _{S(w,\tilde \eta)^c})^{(m)} \} \mathbf 1\{z(\mathbf x) \neq w\}\,\eta^{(m)}(\mathrm{d}\mathbf x)\label{eqn:gineq2}\\
		&\qquad +\frac {1}{m!}\int_{\mathbb{X}^m} g(\mathbf x,\eta) \mathbf 1\{z(\mathbf x) = w\} \,\eta^{(m)}(\mathrm{d}\mathbf x).\label{eqn:gineq3}
	\end{align}
	Now we use that $g$ is stabilizing with stabilization radius $R$. Hence, by \eqref{stoprem} with $\chi:=\eta_{S(w,\tilde \eta)^c}+\tilde \eta_{S(w,\tilde \eta)}$, the second indicator in \eqref{eqn:gineq1} is given by
	\begin{align*}
		\mathbf 1 \{g(\mathbf x,\eta_{S(z(\mathbf x),\eta)}+\eta_{S(z(\mathbf x),\eta)^c \cap S(w,\tilde \eta)^c}+ \tilde \eta_{S(z(\mathbf x),\eta)^c  \cap S(w,\tilde \eta)})\neq g(\mathbf x,\tilde \eta_{S(w,\tilde \eta)}+\eta _{S(w,\tilde \eta)^c})\}. 
	\end{align*}
	Note that if $S(z(\mathbf x),\eta) \cap S(w,\tilde \eta)=\emptyset$, we have $S(z(\mathbf x),\eta)=S(z(\mathbf x),\eta) \cap S(w,\tilde \eta)^c$ and $S(w,\tilde \eta)=S(w,\tilde \eta)\cap S(z(\mathbf x),\eta)^c$, yielding that
	\begin{align*}
		g(\mathbf x,\eta_{S(z(\mathbf x),\eta)}+\eta_{S(z(\mathbf x),\eta)^c \cap S(w,\tilde \eta)^c}+ \tilde \eta_{S(z(\mathbf x),\eta)^c  \cap S(w,\tilde \eta)})= g(\mathbf x,\tilde \eta_{S(w,\tilde \eta)}+\eta _{S(w,\tilde \eta)^c}),
	\end{align*}
	which lets the second indicator in \eqref{eqn:gineq1} vanish. Using Definition  \ref{defstab}(iv) we fnd that for $S(z(\mathbf x),\eta) \cap S(w,\tilde \eta)=\emptyset$ also the first indicators in \eqref{eqn:gineq2} and the indicator in \eqref{eqn:gineq3} vanish. Thus, we conclude that
	\begin{align*}
		(\xi[\eta]-\xi[\tilde \eta_{S(w,\tilde \eta)}+\eta _{S(w,\tilde \eta)^c}]+\delta_w)_+(\mathbb{X})\le	\frac 1 {m!}	\int_{\mathbb{X}^m} g(\mathbf x,\eta) \mathbf 1 \{S(z(\mathbf x),\eta) \cap S(w,\tilde \eta) \neq \emptyset\}\,\eta^{(m)}(\mathrm{d} \mathbf x).
	\end{align*}
	Together with the multivariate Mecke equation \eqref{mecke}, this shows that \eqref{eqn:tv+lhs} is given by
	\begin{align*}
		\frac{1}{(m!)^2}	\int_{\mathbb{X}^m} \int_{\mathbb{X}^m} \mathbb{E} g(\mathbf{x},\eta+\delta_{\mathbf{x}})g(\mathbf{y},\tilde \eta+\delta_{\mathbf{y}})\mathbf 1 \{S(z(\mathbf x),\eta+\delta_{\mathbf x}) \cap S(z(\mathbf y),\tilde \eta+\delta_{\mathbf y}) \neq \emptyset\}\,\lambda^m(\mathrm{d} \mathbf y)\,\lambda^m(\mathrm{d} \mathbf x).
	\end{align*}
	Finally, we want to replace $S(z(\mathbf x),\eta+\delta_{\mathbf x})$ and $S(z(\mathbf y),\tilde \eta+\delta_{\mathbf y})$ by deterministic sets. To achieve this goal, we split he integration area into $\{R(z(\mathbf x),\eta+\delta_{\mathbf x})\le b_{\mathbf{x}},\,R(z(\mathbf y),\tilde \eta+\delta_{\mathbf y})\le b_{\mathbf{y}}\}$ and the complement of this set. Hence, the above is bounded by
	\begin{align*}
		&\frac{1}{(m!)^2}\int_{\mathbb{X}^m} \int_{\mathbb{X}^m} \mathbb{E}g(\mathbf{x},\eta+\delta_{\mathbf{x}})g(\mathbf{y},\tilde \eta+\delta_{\mathbf{y}}) \mathbf 1\{R(z(\mathbf{x}),\eta+\delta_{\mathbf x})>b_{\mathbf{x}} \text{ or } R(z(\mathbf{y}),\tilde\eta+\delta_{\mathbf y})>b_{\mathbf{y}}\} \\
		&\qquad \qquad \qquad \qquad\times \,\lambda^m(\mathrm{d} \mathbf y)\,\lambda^m(\mathrm{d} \mathbf x)\\
		&	\quad +	\frac{1}{(m!)^2}\int_{\mathbb{X}^m} \int_{\mathbb{X}^m} \mathbb{E} g(\mathbf{x},\eta+\delta_{\mathbf{x}})g(\mathbf{y},\tilde \eta+\delta_{\mathbf{y}})\mathbf 1 \{\|z(\mathbf{x})-z(\mathbf{y})\|\le b_{\mathbf{x}}+b_{\mathbf{y}}\}\,\lambda^m(\mathrm{d} \mathbf y)\,\lambda^m(\mathrm{d} \mathbf x)	.
	\end{align*}
	Here, the first term can be bounded by $T_1$ if we assume that $R(z(\mathbf{x}),\eta+\delta_{\mathbf x})>b_{\mathbf{x}}$ (at the cost of a factor 2) and the second term is $T_2$.
\end{proof}

\begin{lemma}\label{lem-}
	Let $\xi=\xi[\eta]$ be the process defined at \eqref{xidef} and for $\mathbb{E} \xi$-almost all $w \in \mathbb{X}$ let $\xi^{w!}$ be the process defined at \eqref{def:xiz}. We have
	\begin{align*}
		\int \mathbb{E}[(\xi-\xi^{w!})_-(\mathbb{X})] \,\mathbb{E}\xi(\mathrm{d}w) \le T_3+T_4+T_5,
	\end{align*}
	where $T_3,T_4,T_5$ are given in Theorem \ref{maithmallg}.
\end{lemma}

\begin{proof}
	Let $\eta$ and $\tilde \eta$ be independent Poisson processes on $\mathbb{X}$ with intensity measure $\lambda$. It follows from the definition of $\xi^{w!}$ in \eqref{def:xiz} that the left-hand side of the statement of the lemma is given by
	\begin{align}
		&\mathbb{E} \int(\xi[\tilde \eta_{S(w,\tilde \eta)}+\eta _{S(w,\tilde \eta)^c}]-\delta_w-\xi[\eta])_+(\mathbb{X})\,\xi[\tilde \eta](\mathrm{d}w)\nonumber\\
		&\quad =\mathbb{E} \int(\xi[\tilde \eta_{S(z(\mathbf x),\tilde \eta)}+\eta _{S(z(\mathbf x),\tilde \eta)^c}]-\delta_{z(\mathbf x)}-\xi[\eta])_+(\mathbb{X})\,g(\mathbf x,\tilde \eta)\,\tilde \eta^{(m)}(\mathrm{d}\mathbf x).\label{eqn:tv-lhs}
	\end{align}
	Now we use that $\lambda$ is diffuse, which implies by \cite[Proposition 6.9]{last2017lectures} that $\eta$ (and hence also $\eta^{(m)}$) is a simple Poisson proecess. Hence, we find
	\begin{align}
		&(\xi[\tilde \eta_{S(z(\mathbf x),\tilde \eta)}+\eta _{S(z(\mathbf x),\tilde \eta)^c}]-\delta_{z(\mathbf x)}-\xi[\eta])_+(\mathbb{X})\nonumber\\
		&\quad \le \frac{1}{m!}\int_{\mathbb{X}^m} g(\mathbf y,\tilde \eta_{S(z(\mathbf x),\tilde \eta)}+\eta _{S(z(\mathbf x),\tilde \eta)^c}) \,\mathbf 1\{\mathbf y \in \eta^{(m)}\}\mathbf 1 \{g(\mathbf{y},\eta) \neq g(\mathbf y,\tilde \eta_{S(z(\mathbf x),\tilde \eta)}+\eta _{S(z(\mathbf x),\tilde \eta)^c})\}\nonumber\\
		&\qquad \qquad \times  \mathbf 1\{\mathbf x \neq \mathbf y\}(\tilde \eta_{S(z(\mathbf x),\tilde \eta)}+\eta _{S(z(\mathbf x),\tilde \eta)^c})(\mathrm{d}\mathbf y)\label{eqn:gineq1-}\\
		&\qquad + \frac{1}{m!}\int_{\mathbb{X}^m} g(\mathbf y,\tilde \eta_{S(z(\mathbf x),\tilde \eta)}+\eta _{S(z(\mathbf x),\tilde \eta)^c}) \,\mathbf 1\{\mathbf y \notin \eta^{(m)}\}  \mathbf 1\{\mathbf x \neq \mathbf y\}(\tilde \eta_{S(z(\mathbf x),\tilde \eta)}+\eta _{S(z(\mathbf x),\tilde \eta)^c})(\mathrm{d}\mathbf y)\label{eqn:gineq2-}.
	\end{align}
	Now we invoke that $g$ is stabilizing and find by \eqref{stoprem} with $\chi:=\eta$ that the second indicator in \eqref{eqn:gineq1-} is given by
	\begin{align*}
		&\mathbf 1 \{g(\mathbf{y},\eta) \neq g(\mathbf y,\tilde \eta_{S(z(\mathbf x),\tilde \eta) \cap S(z(\mathbf y),\omega)}+\eta _{S(z(\mathbf x),\tilde \eta)^c \cap S(z(\mathbf y),\omega)}+\eta _{S(z(\mathbf y),\omega)^c})\},
	\end{align*}
	where $\omega:=\tilde \eta_{S(z(\mathbf x),\tilde \eta)}+ \eta_{S(z(\mathbf x),\tilde \eta)^c}$. Note that the indicator vanishes for $S(z(\mathbf x),\tilde \eta) \cap S(z(\mathbf y),\omega)=\emptyset$. Since in this case also the first indicator in \eqref{eqn:gineq2-} vanishes by Definition \ref{defstab}(iv), we conclude that
	\begin{align*}
		&\mathbb{E} \int(\xi[\tilde \eta_{S(w,\tilde \eta)}+\eta _{S(w,\tilde \eta)^c}]-\delta_w-\xi[\eta])_+(\mathbb{X})\,\xi[\tilde \eta](\mathrm{d}w)\nonumber\\
		&\quad \le \frac{1}{(m!)^2}\mathbb{E}\int_{\mathbb{X}^m}	\int_{\mathbb{X}^m} \mathbf 1 \{S(z(\mathbf x),\tilde \eta) \cap S(z(\mathbf y),\tilde \eta_{S(z(\mathbf x),\tilde \eta)}+\eta_{S(z(\mathbf x),\tilde \eta)^c}) \neq \emptyset\}\,\mathbf 1\{\mathbf x \neq \mathbf y\}\nonumber\\
		&\qquad \qquad \times g(\mathbf{x},\tilde\eta)g(\mathbf y,\tilde \eta_{S(z(\mathbf x),\tilde \eta)}+\eta _{S(z(\mathbf x),\tilde \eta)^c})  \,(\tilde \eta_{S(z(\mathbf x),\tilde \eta)}+\eta _{S(z(\mathbf x),\tilde \eta)^c})^{(m)}(\mathrm{d}\mathbf y) \,\tilde \eta ^{(m)} (\mathrm{d}\mathbf x).
	\end{align*}
	Now we apply the multivariate Mecke equation to the outer integral and obatin
	\begin{align}
		&\frac{1}{(m!)^2}\mathbb{E}\int_{\mathbb{X}^m}	\int_{\mathbb{X}^m}  \mathbf 1 \{S(z(\mathbf x),\tilde \eta+\delta_{\mathbf x}) \cap S(z(\mathbf y),\tilde \eta_{S(z(\mathbf x),\tilde \eta+\delta_{\mathbf x})}+\eta_{S(z(\mathbf x),\tilde \eta+\delta_{\mathbf x})^c}+\delta_{\mathbf x}) \neq \emptyset\}\,\nonumber\\ 
		&\qquad  \times \mathbf 1\{\mathbf x \neq \mathbf y\}\,g(\mathbf x,\tilde \eta+\delta_{\mathbf x})\,g(\mathbf y,\tilde \eta_{S(z(\mathbf x),\tilde \eta+\delta_{\mathbf{x}})}+\eta _{S(z(\mathbf x),\tilde \eta+\delta_{\mathbf{x}})^c}+\delta_{\mathbf{x}})\nonumber\\
		&\qquad \times (\tilde \eta_{S(z(\mathbf x),\tilde \eta+\delta_{\mathbf{x}})}+\eta _{S(z(\mathbf x),\tilde \eta+\delta_{\mathbf{x}})^c}+\delta_{\mathbf{x}})^{(m)}(\mathrm{d} \mathbf y) \lambda^{m} (\mathrm{d}\mathbf x),\label{eqn:tv-meckeout}
	\end{align}
	where we have used that by Definition \ref{defstab}(iv) it holds that $(\tilde \eta+\delta_{\mathbf{x}})_{S(z(\mathbf x),\tilde \eta+\delta_{\mathbf{x}})}=\tilde \eta_{S(z(\mathbf x),\tilde \eta+\delta_{\mathbf{x}})}+\delta_{\mathbf{x}}$.
	Now we use Lemma \ref{stopmecam} with the stopping set $\mu \mapsto S(z(\mathbf x),\mu+\delta_{\mathbf x})$ and with the function $h:\textbf{N}\times \textbf{N} \to [0,\infty)$ that maps $(\mu,\phi)$ to
	\begin{align*}
		\sum\limits_{\mathbf{y} \in (\mu+\phi+\delta_\mathbf{x})^{(m)}}\mathbf{1}\{S(z(\mathbf{x}),\mu+\delta_\mathbf{x})\cap S(z(\mathbf{y}),\mu+\phi+\delta_\mathbf{x})\neq \emptyset\} \mathds{1}\{\mathbf x \neq \mathbf{y}\}\, g(\mathbf{x},\mu)\,g(\mathbf{y},\mu+\phi+\delta_\mathbf{x}),
	\end{align*}
	we find that \eqref{eqn:tv-meckeout} is given by
	\begin{align*}
		&\frac{1}{(m!)^2}\mathbb{E}\int_{\mathbb{X}^m}	\int_{\mathbb{X}^m}  \mathbf 1 \{S(z(\mathbf x), \eta+\delta_{\mathbf{x}})\cap S(z(\mathbf y),\eta+\delta_{\mathbf{x}}) \neq \emptyset\}\,\mathbf 1\{\mathbf x \neq \mathbf y\} g(\mathbf{x},\eta+\delta_{\mathbf{x}})g(\mathbf y,\eta+\delta_{\mathbf{x}})\nonumber\\
		&\qquad  \qquad\times (\eta+\delta_{\mathbf{x}})^{(m)}(\mathrm{d} \mathbf y) \lambda^{m} (\mathrm{d}\mathbf x).
	\end{align*}
	Here we distinguish by the number $0\le\ell \le m-1$ of elements that $\mathbf{x}$ and $\mathbf{y}$ have in common (note that $\ell=m$ is not possible since $\mathbf{x}\neq \mathbf{y}$). This gives for the above
	\begin{align*}
		&\sum\limits_{\ell=0}^{m-1}\frac{1}{m!\,(m-\ell)!}\mathbb{E}\int_{\mathbb{X}^m}	\int_{\mathbb{X}^{m-\ell}} \mathbf 1 \{S(z(\mathbf x_\ell,\mathbf y), \eta+\delta_{\mathbf{x}})\cap S(z(\mathbf y),\eta+\delta_{\mathbf{x}}) \neq \emptyset\} \,\nonumber\\
		&\qquad \times g((\mathbf{x}_\ell,\mathbf{y}),\eta+\delta_{\mathbf{x}})\,g(\mathbf y,\eta+\delta_{\mathbf{x}})\,\eta^{(m-\ell)}(\mathrm{d}\mathbf y) \lambda^{m} (\mathrm{d}\mathbf x)\nonumber\\
		&\quad=\sum\limits_{\ell=0}^{m-1}\frac{1}{m!\,(m-\ell)!}\mathbb{E}\int_{\mathbb{X}^m}	\int_{\mathbb{X}^{m-\ell}}  \mathbf 1 \{S(z(\mathbf x_\ell,\mathbf y), \eta+\delta_{(\mathbf{x},\mathbf y})\cap S(z(\mathbf y),\eta+\delta_{(\mathbf{x},\mathbf y)}) \neq \emptyset\} \,\nonumber\\
		&\qquad \times g((\mathbf{x}_\ell,\mathbf{y}),\eta+\delta_{(\mathbf{x},\mathbf y)})\,g(\mathbf y,\eta+\delta_{(\mathbf{x},\mathbf y)})\lambda^{m-\ell}(\mathrm{d}\mathbf y) \lambda^{m} (\mathrm{d}\mathbf x).
	\end{align*}
	Here, all terms with $\ell \neq 0$ form the term $T_5$ from the statement of the lemma. For $\ell =0$ we distinguish by the sizes of the stabilization radii. This gives the terms $T_3$ and $T_4$.
\end{proof}

\begin{proof} [Proof of Theorem \ref{maithmallg}] 
	By Lemma \ref{lemxiz}, the process $\xi^{w!}$ defined at \eqref{def:xiz} is for $\mathbb{E}\xi$-almost all $w \in \mathbb{X}$ a Palm version of $\xi$. Hence, we find from \cite[Theorem 2.6]{barbour1992stein} that
	\begin{align*}
		\mathbf{d_{TV}} (\xi, \,\nu) \le \|\mathbb{E} \xi- \mathbb{E} \nu \| + \int\limits \E \|\xi-\xi^{w!}\|\,\mathbb E \xi(\mathrm{d}w).
	\end{align*}
	Now we invoke Lemma \ref{lem+} and Lemma \ref{lem-} to bound the integral. This finishes the proof of Theorem \ref{maithmallg}.
\end{proof}

\section{Maximum cells in the Poisson-Voronoi mosaic} \label{s5}
In this section we apply Theorem \ref{maithmallg} to point processes of centres of large cells in the Poisson-Voronoi mosaic. Let $\mathbb{X}=\mathbb{R}^d$ ($d \ge 2$) with Borel $\sigma$-field $\mathcal{B}^d$ with $d$-dimensional Lebesgue measure $\L_d$ and standard scalar product $\langle \cdot,\cdot \rangle$. For $\mu \in \mathbf{N}_s$ and $x \in \mu$ the {\em Voronoi cell} $C(x,\mu)$ is the set of all points $y \in \R^d$ with $|y-x|\le \min_{z \in \mu} |y-z|$, where $|\cdot|$ is the Euclidean norm. It is a closed convex set with interior points. If $\mu=\eta$ is a (stationary) Poisson process in $\mathbb{R}^d$, the system $\{C(x,\mu):\,x \in \eta\}$ is called {\em Poisson-Voronoi mosaic}. For an in-depth introduction to the theory of (random) mosaics we refer to Section 10 in \cite{schneider2008stochastic}. 

Next we explain how we measure the size of a Voronoi cell. Let $\mathcal{K}_o^d$ denote the space of all convex bodies (nonempty, compact and convex sets) $K \subset \mathbb{R}^d$ containing the origin $o$ as an interior point and equip $\mathcal{K}_o^d$ with the Hausdorff metric. Following \cite{hug2007asymptotic} let $k>0$ and call a map  $\Sigma: \mathcal{K}_o^d \to [0,\infty)$ {\em size functional} if $\Sigma$ is continuous, not identically $0$, $k$-homogeneous (i.e.\ $\Sigma(aK)=a^k \Sigma(K)$ for all $a>0$ and $K \in \mathcal{K}_o^d$, where $aK:=\{ax:\,x \in K\}$) and increasing under set inclusions (i.e.\ $\Sigma(K_1) \le \Sigma(K_2)$ for all $K_1,K_2 \in \mathcal{K}_o^d$ with $K_1 \subset K_2$). The {\em centred inradius} $\rho_o$ (i.e.\ the inradius of the largest ball with centre $o$ contained in $K$) is an example for a $1$-homogeneous size functional. The {\em $k$th intrinsic volume} $V_k$ ($k \in [d]$) of $K$ serves as an example for a $k$-homogeneous size functional. These size functionals are discussed in more details in Example \ref{PVTEx1}. 

We study cells in the Poisson-Voronoi mosaic that are large with respect to a size functional $\Sigma$.
Let $c>0$, $W \subset \mathbb{R}^d$ be compact, $\eta_\gamma$ be a stationary Poisson process with intensity $\gamma>0$ and let $\Sigma$ be a $k$-homogeneous size functional. We slightly abuse the notation and write $\Sigma(C(x,\mu)):=\Sigma(C(x,\mu)-x)$ for $\mu \in \mathbf N_s$ with $x \in \mu$. For a threshold $v_{c,\gamma}>0,\,\gamma > 0,$ (to be specified in {\eqref{PVT:defvn} below) we consider the process
	\begin{align}
		\xi_{c,\gamma}:=\sum\limits_{x \in \eta} \mathbf{1}\{\Sigma(C(x,\eta_\gamma))>v_{c,\gamma}\}\,\delta_{x}. \label{PVTxi}
	\end{align}
	Here, the threshold $v_{c,\gamma}$ is chosen such hat the intensity measure $\mathbb{E}\xi_{c,\gamma}$ of $\xi_{c,\gamma}$ satisfies
	\begin{align}
		\mathbb{E} \xi_{c,\gamma}(A)=c\lambda_d(A),\quad\gamma >0,\, A \in \mathcal B^d. \label{PVTintvt}
	\end{align}
	Note that by the Mecke equation \eqref{mecke} and by stationarity of $\eta$ we have 
	\begin{align}
		\mathbb{E} \xi_{c,\gamma}(A)=\gamma \lambda_d(A)\,\P(\Sigma(C(o,\eta_\gamma+\delta_o))>v_{c,\gamma}).\label{PVTmecke}
	\end{align}
	To see that $v_{c,\gamma}$ can be chosen such that \eqref{PVTxi} exists, note that by \cite[Section 9]{hug2007asymptotic} the distribution $\mathbb{P}^{\Sigma(C(o,\eta_\gamma+\delta_o))}$ of $\Sigma(C(o,\eta_\gamma+\delta_o))$ and the Lebesgue meaure $\lambda_1$ are equivalent measures on $[0,\infty)$. Together with \eqref{PVTmecke}, this implies that the choice
	\begin{align}
		v_{c,\gamma}:=\inf\{v>0:\,\gamma\mathbb{P}(\Sigma(C(o,\eta_\gamma+\delta_o))>v)>c\},\quad \gamma >0,\label{PVT:defvn}
	\end{align}
	indeed satisfies \eqref{PVTintvt}. 
	
	We need to introduce some more notation. For $K \in \mathcal{K}_o^d$ let $h_K(u):=\max\{\langle x,u \rangle:\,x \in K\},\,u\in \mathbb S^{d-1},$ be the {\em support function} of $K$. Define
	\begin{align*}
		\Phi(K):=\frac 1d \int_{\mathbb S^{d-1}}h_K(u)^d\, \sigma(\mathrm{d}u),
	\end{align*}
	where we write $\sigma$ for the uniform distribution on the unit sphere $\mathbb{S}^{d-1}$ in $\R^d$. There is a constant $\tau >0$ such that $\Phi$ and $\Sigma$ satisfy the sharp isoperimetric inequality
	\begin{align} \label{lowbvor}
		\Phi(K)\ge \tau \Sigma(K)^{d/k},\quad K \in \mathcal{K}_o^d.
	\end{align}
	That this inequality is sharp means that there is some $K \in \mathcal{K}_o^d$ with more that one point for which equality holds in \eqref{lowbvor} (see \cite[Section 3]{hug2007asymptotic}). Every such body is called an {\em  extremal body}. For example, if $\Sigma$ is the $d$-dimensional volume, $\tau=(d \kappa_d)^{-1}$ and the extremal bodies are exactly the $d$-dimensional balls centred at the origin $o$.
	
	We call a non-negative, continuous, $0$-homogeneous functional $\vartheta: \mathcal{K}_o^d \to [0,\infty)$ {\em deviation functional} if it has the property that $\vt(K)=0$ holds for $K \in \mathcal{K}_o^d$ with $\Sigma(K)>0$ if and only if $K$ is an extremal body. Such deviation functionals always exist. For example, 
	\begin{align*}
		\vartheta(K):=\frac{\Phi(K)}{\tau \Sigma(K)^{d/k}}-1
	\end{align*}
	from (5) in \cite{hug2007asymptotic} defines a deviation functional. There exists a continuous function $f:(0,+\infty) \to (0,+\infty)$ with $f(0)=0$ and $f(\eps)>0$ for $\eps>0$ such that
	\begin{align*}
		\Phi(K)\ge (1+f(\eps)) \tau \Sigma(K)^{d/k} \quad \text{for } \vt(K)\ge \eps,
	\end{align*}
	which sharpens \eqref{lowbvor}. Any such function is called a {\em stability function}.
	
	The following statement is Theorem 1 in \cite{hug2007asymptotic} (specialized to the Poisson-Voronoi mosaic). Suppose that a stationary Poisson process $\eta$ with intensity $\gamma$, a size functional $\Sigma$, a deviation functional $\vt$ and a stability function $f$ (for $\Phi$, $\Sigma$ and $\vt$) are given. Then the following holds. There exists a positive constant $c_0$ (depending only on $\tau$) such that for all $\eps>0$ and $v>0$ we have
	\begin{align}
		\mathbb{P}(\vt(C(o,\eta_\gamma+\delta_o)) \ge \eps \mid \,\Sigma(C(o,\eta_\gamma+\delta_o))>v) \le c_1 \exp(-c_0f(\eps) v^{d/k}\gamma),\label{PVT:shape}
	\end{align}
	where $c_1>0$ depends only on $d$, $\Sigma$, $f$, $\eps$.
	
	Next we determine the asymptotic behavior of $v_{c,\gamma}$ as $\gamma \to \infty$. Since $v_{c,\gamma}\to \infty$ as $\gamma \to \infty$ (which follows from the definitoon of $v_{c,\gamma}$ together with the equivalence of $\mathbb{P}^{\Sigma(C(o,\eta_\gamma+\delta_o))}$ and $\lambda_1$ on $[0,\infty)$) we find from Theorem 2 in \cite{hug2007asymptotic} that
	\begin{align*}
		\lim_{\gamma \to \infty} \gamma^{-1}v_{c,\gamma}^{-d/k}\log \P(\Sigma(C(o,\eta_\gamma+\delta_o))>v_{c,\gamma}) =-2^d d \kappa_d \tau,\quad c>0,
	\end{align*}
	where $\tau$ is the constant from \eqref{lowbvor}. Since $\gamma^{-1}v_{c,\gamma}^{-d/k}\log \gamma$ is given by
	\begin{align*}
		\gamma^{-1} v_{c,\gamma}^{-d/k}\log \P(\Sigma(o,\eta_\gamma+\delta_o)>v_{c,\gamma}) \Big(\frac{\log(\gamma\P(\Sigma(o,\eta_\gamma+\delta_o)>v_{c,\gamma}))}{\log \gamma}-1\Big)^{-1}
	\end{align*}
	and since $\gamma\P(\Sigma(o,\eta_\gamma+\delta_o)>v_{c,\gamma})=c$ for all $\gamma>0$ we conclude that $\gamma^{-1}v_{c,\gamma}^{-d/k}\log \gamma\to 2^d d\kappa_d \tau$ as $\gamma \to \infty$.
	
	In our Poisson approximation result for $\xi_{c,\gamma}$ we will need the following condition on the extremal bodies of the size functional $\Sigma$. For $K \in \mathcal{K}_o^d$ let $r_o(K)$ be the radius of the smallest ball with center $o$ containing $K$ (centred circumradius) and $\rho_o(K)$ be the radius of the largest ball with center $o$ contained in $K$ (centred inradius). We assume that there exists a function $h:[0,\infty)\to (0,1]$ such that for some $\varepsilon >0$:
	\begin{align} 
		\frac{\rho_o(K)}{r_o(K)} \ge h(\eps) \quad\text{for all }K\in \mathcal{K}_o^d \text{ with }\vt(K)<\eps.\label{PVTass1}
	\end{align}
	
	The following theorem is the main result of this section.
	
	\begin{theorem} \label{PVTmaith}
		Suppose that a stationary Poisson process $\eta_\gamma$ with intensity $\gamma>0$, a size functional $\Sigma$, a deviation functional $\vt$  and a stability function $f$ are given and that \eqref{PVTass1} holds for some $\eps>0$. Let $W \subset \mathbb{R}^d$ be compact, $c>0$ and let $\nu_c$ be a stationary Poisson process with intensity $c$. Then we have for all $\delta>0$
		\begin{align*}
			\mathbf{d_{TV}}(\xi_{c,\gamma} \cap W,\,\nu_c\cap W) \le C \gamma^{\delta-\min\big[c_0 f(\eps), h(\eps)^d  k(h(\eps))\big]},\quad \gamma>0,
		\end{align*}
		where the constant $C>0$ does not depend on $\gamma$. Here, $c_0$ is the constant from \eqref{PVT:shape} and $k(a)$ is the volume of the intersection of an infinite cone with apex $o$ and angular radius $\arcsin\big(\frac{a}{\sqrt{1+a^2}}\big)$ and of a $d$-dimensional ball centred at $o$ with volume 1.
	\end{theorem}
	
	As the proof will show, the expoenent of $\gamma$ of the right-hand side of the statement in Theorem \ref{PVTmaith} has a clear geometric interpretation. While the first term $c_0f(\eps)$ comes from \eqref{PVT:shape} and can be interpreted as the approximation error of a large Voronoi cell by an extremal body, the second term $h(\eps)^d  k(h(\eps))$ stems from a stabilization result for Voronoi cells whose shape is close to that of an extremal body (see Lemma \ref{stabaug}).
	
	We now demonstrate how Theorem \ref{PVTmaith} applies to concrete size functionals $\Sigma$ and how it can be used to derive extreme value statements for large cells in the Poisson-Voronoi mosaic.
	
	\begin{example} \label{PVTEx1} (a) Let $\Sigma:=\rho_o$ be the centred inradius. Then $\tau=1/d$ and the isoperimetric inequality \eqref{lowbvor} reads $\rho_o(K)^d\le d\Phi(K),\, K \in \mathcal{K}_o^d,$ where equality holds if and only if $K$ is a $d$-dimensional ball centred at the origin $o$. We choose the deviation functional
		\begin{align*}
			\vt(K):=\frac{r_o(K)-\rho_o(K)}{r_o(K)+\rho_o(K)} \in [0,1].
		\end{align*}
		Hence, for $\eps \in (0,1)$ we have that $\vt(K)< \eps$ if and only if $\frac{\rho_o(K)}{r_o(K)}> \frac{1-\eps}{1+\eps}$. Therefore, \eqref{PVTass1} holds with $h(\eps):=\frac{1-\eps}{1+\eps}$ and a stability function is given by $f(\eps):=\big(\frac{1+\eps}{1-\eps}\big)^d-1$. Since
		\begin{align*}
			\{\rho_o(C(x,\eta_\gamma+\delta_x))>v\}=\{\eta_\gamma \cap B_{2v}(x)=\varnothing\}
		\end{align*}
		we choose $v_{c,\gamma}^d:=2^{-d}\kappa_d^{-1} \gamma^{-1} \log (\gamma/c)$. Letting $c:=e^{-t}$ for $t \in \mathbb{R}$ and $$M:=\min\Big[c_0\Big(\big(\frac{1+\eps}{1-\eps}\big)^d-1\Big), \big(\frac{1-\eps}{1+\eps}\big)^d k\Big(\frac{1-\eps}{1+\eps}\Big)\Big]>0$$ we find from Theorem \ref{PVTmaith} that for all $\delta>0$
		\begin{align*}
			\Big|\mathbb{P}(2^d \kappa_d \gamma \max_{x \in \eta_\gamma \cap W}\rho_o(C(x,\eta_\gamma))^d-\log \gamma  \le t )-\exp(-\lambda_d(W)e^{-t})\Big|\le C \gamma^{\delta-M},\quad t \in \mathbb{R},\,\gamma>0.
		\end{align*} 
		This shows that $\max_{x \in \eta_\gamma \cap W}\rho_o(C(x,\eta))^d$ is in the domain of attraction of Gumbel distribution and quantifies the rate of convergence in  statement (2a) from Theorem 1 \cite{calka2014extreme}.\\
		(b) Let $\Sigma:=V_k$ ($1 \le k \le d$) be the $k$th intrinsic volume (in particular, $V_d$ is the volume, $dV_{d-1}$ is the surface area and $2V_1/\kappa_d$ is the mean width). From (15) in \cite{hug2007asymptotic} we have that
		\begin{align*}
			\frac 1d \Big(\frac{ k! (d-k)!\kappa_{d-k}}{ d!\kappa_d} \Big)^{d/k} V_k(K)^{d/k} \le \Phi(K),\quad K \in \mathcal{K}_o^d.
		\end{align*}
		As in (a), the extremal bodies are precisely the $d$-dimensional balls with centre at $o$. Hence, $\vt$ and $a$ can be chosen as above.
	\end{example}
	
	As a preparation for the proof of Theorem \ref{maithmallg} we show that the  function
	\begin{align}
		g(x,\mu)=\mathbf 1\{x \in W\} \mathbf 1\{\Sigma(C(x,\mu))>v_{c,\gamma}\},\quad x \in \mathbb{R}^d,\,\mu \in \mathbf{N}_s, \label{gvor}
	\end{align}
	is stabilizing in the sense of Definition \ref{defstab} and we construct a stabilization radius that satisfies the conditions from Definiton \ref{defstab}. Following Section 6.3 in \cite{penrose2007gaussian} let $K_i(x),\,1 \le i \le I,$ be a finite collection of infinite open cones in $\R^d$ with angular radius $\pi/6$, apex at $x$ and union $\R^d$. For $\mu \in \mathbf {N}_s$ we define 
	\begin{align}
		R_i(x,\mu):&=\inf\{r>0:\,K_i(x)\cap B_r(x)\neq\emptyset\},\quad i \in I,\label{stabivor}\\
		R(x,\mu):&=2\max\limits_{1 \le i \le I} R_i(x,\mu).\label{stabvor}
	\end{align}
	Then we have 
	\begin{align*}
		C(x,\mu+\delta_x)=C(x,\mu_{B(x,2R(x,\mu))}+\delta_x)
	\end{align*}
	which implies that $g$ from \eqref{gvor} is stabilizing. Since $\P(R(x,\eta_\gamma)<\infty)=1$ and since $\mu \mapsto B(x,R(x,\mu))$ is a stopping set we conclude that $R$ is a stabilization radius. 
	By construction we find that
	\begin{align}
		\P(R(o,\eta_\gamma+\delta_o)>r)\le 1-(1-e^{-\gamma (r/2)^d/I})^{I}\sim e^{-\gamma (r/2)^d/I} \quad \text{as }r\to \infty.\label{PVTRass}
	\end{align}

	Next we derive a more refined stabilization property that holds in a mosaic for which is underlying point configuration is augmented by a given element $y \in \mathbb{R}^d$. Let $\mathbb S_{\{y-x\}^\perp}$ denote the unit sphere in the linear subspace orthogonal to $y-x$. We define the infinite cone
	\begin{align*}
		K_a(x,y):=\{y+t(y-x) +at |y-x|u:\,t> 0,\,u \in \mathbb S_{\{y-x\}^\perp} \}.
	\end{align*}
	Hence, $K_a(x,y)$ has apex $o$, axis $y-x$ and angular radius $\arcsin\big(\frac{a}{\sqrt{1+a^2}}\big)$.
	
	In the proof of our main theorem of this section we will make use of the following statement.
	\begin{lemma}\label{stabaug}
		Let $x,y \in \mathbb{R}^d$ such that $\mu+\delta_{(x,y)} \in \mathbf{N}_s$. Let $\vt(C(x,\mu+\delta_{(x,y)}))\le \eps$ for some $\eps>0$ and assume that condition \eqref{PVTass1} holds for $\eps>0$. Then we have
		\begin{align*}
			C(x,\mu+\delta_{(x,y)})=C(x,\mu \cap K_a(x,y)^c+\omega \cap K_a(x,y)+\delta_{(x,y)}),\quad \omega \in \mathbf N_s,
		\end{align*}
		where $a:=\frac{h(\eps)}{\sqrt{1-h(\eps)^2}}$ with $h$ from \eqref{PVTass1}.
	\end{lemma}
	
	\begin{figure}[bt]
		\centering
		\includegraphics[scale=0.39]{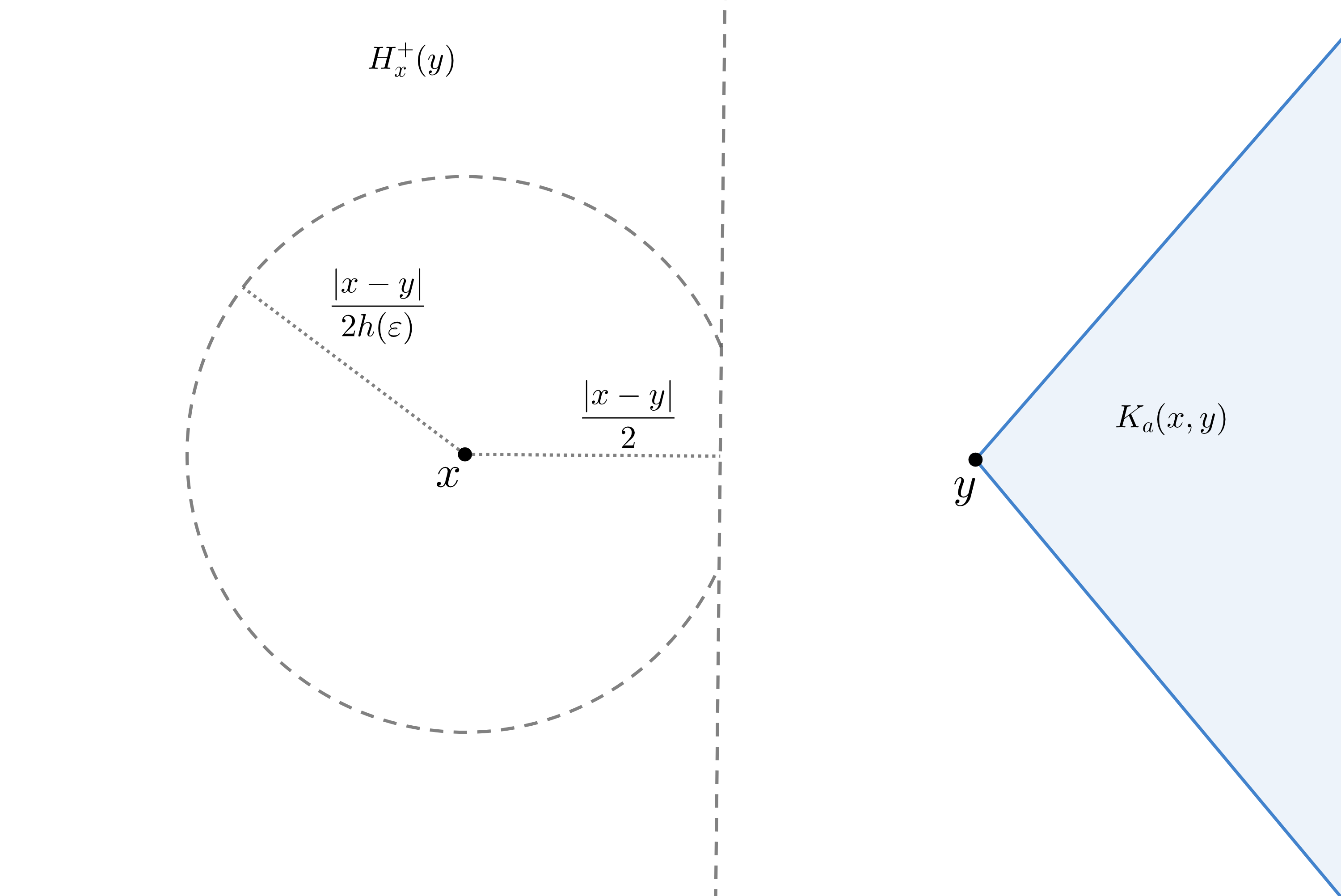}
		\caption{The cone $K_a(x,y)$ is marked in blue.}
	\end{figure}
	
	\begin{proof}
		First note that $\rho_o(C(x,\mu+\delta_{(x,y)}))\le \frac{|x-y|}{2}$. By \eqref{PVTass1}, this together with the condition $\vt(C(x,\mu+\delta_{(x,y)}))\le \eps$ yield that $r_o(C(x,\mu+\delta_{(x,y)}))\le \frac{|x-y|}{2h(\eps)}$. Hence, it suffices to show that
		\begin{align}
			|z-w| \ge \frac{|x-y|}{2h(\eps)}\quad \text{for all } z \in \{y\} \cup K_a(x,y),\,w\in \partial B\Big(x,\frac{|x-y|}{2h(\eps)}\Big) \cap H_x^+(x),\label{stabaughalv}
		\end{align}
		where $H_x^+(y)$ is the closed half-space that contains $x$ and which is delimited by the bisecting hyperplane of $[x,y]$. We use the parametrizations $z=y+t_1(y-x)+t_2|x-y|v$ and $w=y+s_1(y-x)+s_2|x-y|u$ with 
		$$s_1 \in \Big[-1-\frac{1}{2h(\eps)},-\frac 12\Big],\quad s_2=\sqrt{\frac{1}{4h(\eps)^2}-(s_1+1)^2}, \quad t_1\ge 0,\quad t_2 \in \Bigg[0,\frac{h(\eps)}{\sqrt{1-h(\eps)^2}}t_1\Bigg]$$
		and $u,v \in S_{\{y-x\}^\perp}$. We have
		\begin{align*}
			\frac{|z-w|^2}{|x-y|^2} &\ge (t_1-s_1)^2+(t_2-s_2)^2
		\end{align*}
		which attains its minimum value $\frac{1}{4h(\eps)^2}$ for $t_1=0$, $t_2=0$ and $s_1=-\frac 12$ (under the constraints on $s_1$, $s_2$, $t_1$ and $t_2$ from above). This shows \eqref{stabaughalv}.
	\end{proof}

	\begin{proof}[Proof of Theorem \ref{PVTmaith}]
		For $\gamma>0$ we apply Theorem \ref{maithmallg} with the function $g$ from \eqref{gvor}, $z(x):=x$, $\lambda:=\gamma \lambda_d \cap W$ and the stabilization radius $R$ defined at \eqref{stabvor}. Let
		\begin{align}
			b_x:=b_\gamma=2(3I \gamma^{-1} \log \gamma)^{1/d},\quad x \in \mathbb{R}^d, \label{PVTdefbt}
		\end{align}
		and let $\nu_c$ be a stationary Poisson process in $\mathbb{R}^d$ with intensity $c$. From Theorem \ref{maithmallg} (with $T_5=0$ since $m=1$) we obtain that
		\begin{align*}
			\mathbf{d_{TV}}(\xi_{c,\gamma} \cap W,\nu_{c,\gamma} \cap W)\le \|\mathbb{E} (\xi_{c,\gamma}\cap W)-\mathbb{E}(\nu_{c}\cap W)\|+T_1+T_2+T_3+T_4
		\end{align*}	
		with
		\begin{align*}
			T_1&=2\gamma \mathbb{E} \xi_{c,\gamma}(W) \int_{W} \mathbb{P}(\Sigma(C(x,\eta_\gamma+\delta_x))>v_{c,\gamma},\,R(x,\eta_\gamma+\delta_x)>b_\gamma)\,\mathrm{d}x,\\
			T_2&=\gamma^2\int_{W} \int_{W}\,\mathbf 1\{|x-y|\le 2b_\gamma\}\mathbb{P}(\Sigma(C(x,\eta_\gamma+\delta_x))>v_{\gamma}) \mathbb{P}(\Sigma(C(y,\eta_\gamma+\delta_y))>v_{n})\,\mathrm{d}y\,\mathrm{d}x,\\
			T_3&=\gamma^2\int_{W} \int_{W} \mathbb{P}(\Sigma(C(x,\eta_\gamma+\delta_x))>v_{\gamma},\,\Sigma(C(y,\eta_\gamma+\delta_y))>v_{c,\gamma},\,R(x,\eta_\gamma+\delta_x)>b_\gamma)\,\mathrm{d}y\,\mathrm{d}x,\\
			T_4&=\gamma^2\int_{W} \int_{W}\mathbf 1\{|x-y|\le 2b_\gamma\}\, \mathbb{E} \mathbf 1\{\Sigma(C(x,\eta_\gamma+\delta_{(x,y)}))>v_{c,\gamma},\,\Sigma(C(y,\eta_\gamma+\delta_{(x,y)}))>v_{c,\gamma}\}\\
			&\qquad \times \mathbf 1\{R(x,\eta_\gamma+\delta_{(x,y)})\le b_\gamma,\,R(y,\eta_\gamma+\delta_{(x,y)})\le b_\gamma\}\,\mathrm{d}y\,\mathrm{d}x.
		\end{align*}
		where $\mathrm d :=\mathrm{d}\lambda_d$ denotes integration with respect to the Lebesgue measure $\lambda_d$. In the following, $\alpha_i>0$ ($i \in \mathbb{N}$) are positive constants that do not depend on $\gamma$. Their precise values are not important for the argument.

		First we show that the total variation on the right-hand side vanishes. Note that by the Mecke equation, by stationarity of $\eta$ and by \eqref{PVTintvt} we have 
		\begin{align*}
			\mathbb{E} \xi_{c,\gamma}(A)=n \lambda_d(A) \mathbb{P}(\Sigma(C(o,\eta_\gamma+\delta_o))>v_{c,\gamma})=c\lambda_d(A)=\nu_c(A),\quad A \in\mathcal B^d,
		\end{align*}
		which yields that the intensity measures of $\xi_{c,\gamma}$ and $\nu_{c}$ coincide for all $\gamma>0$ and $c>0$.
		
		{\em The estimate of $T_1$.} We bound the probability in the integral of $T_1$ by
		\begin{align*}
			\P(R(x,\eta_\gamma+\delta_x)>v_{c,\gamma})\le 1-\prod_{i \in I} \mathbb{P}(R_i(x,\eta_\gamma+\delta_x)>b_\gamma)=1-(1-\mathrm e^{-\gamma (b_\gamma/2)^d/I})^{I}\le Ie^{-\gamma (b_\gamma/2)^d/I}
		\end{align*}
		where $R_i(x,\mu)$ is defined at \eqref{stabivor}. Using the definition of $b_\gamma$, we find that
		\begin{align*}
			T_1 \le \frac{2 \gamma\mathbb{E} \nu_{c}(W)\lambda_d(W)}{\gamma^3}\le \frac{\alpha_1}{\gamma^2}.
		\end{align*}
		
		{\em The estimate of $T_2$.} Since $\Sigma$ is translation-invariant and $\eta_\gamma$ is stationary, $T_2$ is bounded by
		\begin{align*}
			T_2 \le \gamma^2 \lambda_d(W) \lambda_d(B_{2b_\gamma}(o)) \mathbb{P}(\Sigma(C(o,\eta_\gamma+\delta_o))>v_{c,\gamma})^2\le \frac{\alpha_2 \log \gamma}{\gamma}
		\end{align*}
		where we have used the definitions of $b_\gamma$ and of $v_{c,\gamma}$ to obtain the second inequality.
		
		{\em The estimate of $T_3$.} We estimate the probability in the integral of $T_3$ in the same way as we did for $T_1$. This yields the bound
		\begin{align*}
			T_3 \le \gamma^2 \lambda_d(W)^2 e^{-\gamma (v_{c,\gamma}/2)^d/I}\le \frac{\alpha_3}{\gamma}.
		\end{align*}
		
		{\em The estimate of $T_4$.} Let $\eps>0$ be the constant from Theorem \ref{PVTmaith}. We distinguish by the shape of $C(x,\eta_\gamma+\delta_{(x,y)})$ which gives the bound
		\begin{align}
			&2\gamma^2 \int_{W} \int_{W} \mathbf 1\{|x-y|\le 2b_\gamma\} \mathbb{P}(\Sigma(C(x,\eta_\gamma+\delta_{(x,y)})>v_{c,\gamma},\,\vt(C(x,\eta_\gamma+\delta_{(x,y)})>\varepsilon)\,\mathrm{d}y\,\mathrm{d}x\label{eqn:vorT4a}\\
			&\quad +  \gamma^2\int_{W} \int_{W} \mathbf 1\{|x-y|\le 2b_\gamma\} \, \mathbb{E}\,\mathbf 1\{\Sigma(C(x,\eta_\gamma+\delta_{(x,y)}))>v_{c,\gamma},\,\vt(C(x,\eta_\gamma+\delta_{(x,y)}))<\varepsilon\}\nonumber\\
			&\qquad \qquad \times \mathbf 1\{\Sigma(C(y,\eta_\gamma+\delta_{(x,y)}))>v_{c,\gamma}\,\vt(C(y,\eta_\gamma+\delta_{(x,y)}))<\varepsilon\}\,\mathrm{d}y\,\mathrm{d}x.\label{eqn:vorT4b}
		\end{align}
		
		Since the probability in \eqref{eqn:vorT4a} uses the augmented process $\eta+\delta_{(x,y)}$ instead of $\eta+\delta_y$, we can not directly invoke \eqref{PVT:shape}. Instead, we apply the Mecke formula to the inner integral and use that $\eta_\gamma$ is stationary. This gives
		\begin{align*}
			2\gamma^2 \lambda_d(W) \mathbb{E}\sum_{y \in \eta \cap B_{2\gamma}} \mathbf 1\{C(o,\eta_\gamma+\delta_{o})>v_{c,\gamma},\,\vt(C(o,\eta_\gamma+\delta_{o}))>\varepsilon\}
		\end{align*}
		Next distinguish by the number $B_\gamma:=\eta (B_{2b_\gamma})$ of points of $\eta_\gamma$ in $B_{2b_\gamma}$. This yields \eqref{eqn:vorT4a} the bound
		\begin{align}
			&2\gamma^2 \lambda_d(W) \,\mathbb{E} B_\gamma\mathbf 1\{B_\gamma \le 2\mathbb{E} B_\gamma\}  \mathbf1 \{\Sigma(C(o,\eta_\gamma+\delta_{o}))>v_{c,\gamma},\,\vt(C(o,\eta_\gamma+\delta_{o}))>\varepsilon\}\label{eqn:vorT4aa}\\
			&\quad +2\gamma^2 \lambda_d(W) \,\mathbb{E} B_\gamma\mathbf 1\{B_\gamma > 2\mathbb{E} B_\gamma\}  \mathbf1 \{\Sigma(C(o,\eta_\gamma+\delta_{o}))>v_{c,\gamma},\,\vt(C(o,\eta_\gamma+\delta_{o}))>\varepsilon\}.\label{eqn:vorT4ab}
		\end{align}
		Since $\mathbb{E} B_{\lambda}=\gamma(2b_\gamma)^d \kappa_d =4^d 3I \log (\gamma) \kappa_d$ we obtain for \eqref{eqn:vorT4aa} the bound
		\begin{align*}
			&4\gamma^2 \lambda_d(W_n) \,\mathbb{E} B_\gamma\mathbb{P}(\Sigma(C(o,\eta_\gamma+\delta_{o}))>v_{c,\lambda})\,\mathbb{P}(\vt(C(o,\eta_\gamma+\delta_{c,\gamma}))>\varepsilon \mid \,\Sigma(C(o,\eta_\gamma+\delta_{o}))>v_{c,\gamma})\\
			&\quad \le \alpha_4 \log(\gamma) \exp(-c_0 f(\varepsilon) v_{c,\gamma}^{d/k}\gamma),
		\end{align*}
		where we have used \eqref{PVT:shape} to obtain the inequality. For \eqref{eqn:vorT4ab} we find by Cauchy-Schwarz the bound
		\begin{align}
			2 \gamma^2 \lambda_d(W) \sqrt{\mathbb{E}B_\gamma} \sqrt{\mathbb P(B\gamma >2 \mathbb{E}B_\gamma)}. \label{eqn:vorT4abcher}
		\end{align}
		Recall that $B_\gamma$ follows a Poisson distribution with parameter $\mathbb{E} B_\gamma$. Hence, we find from the Chernoff bound \cite[Section 5.3]{mitzenmacher2017probability}
		\begin{align}
			\P(X>u)\le \frac{(\mathrm e \lambda)^u \mathrm e^{-\lambda}}{\mathrm e^{u \log u}},\quad u >\lambda,\label{cher}
		\end{align}
		where $X$ is Poisson distributed with parameter $\lambda>0$, that \eqref{eqn:vorT4abcher} is bounded by
		\begin{align*}
			2 \gamma^2 \lambda_d(W) \sqrt{\mathbb{E}B_\gamma} \sqrt{\frac{(\mathrm e \mathbb{E} B_\gamma)^{2 \mathbb{E}B_\gamma}\mathrm e^{-\mathbb{E}B_\gamma}}{\mathrm e ^{2 \mathbb{E} B_\gamma \log (2 \mathbb{E} B_\gamma)}}} &\le 2 \gamma^2\lambda_d(W) \sqrt{\mathbb{E}B_\gamma} \mathrm e^{-\mathbb E B_\gamma (1/2-\log 2)}\\
			&\le \alpha_5  \gamma^{2-4^d3I\kappa_d(2-\log 2)}\sqrt{\log \gamma},
		\end{align*}
		where we have used that $\mathbb{E} B_{\lambda}=4^d 3I \log (\gamma) \kappa_d$.
		
		To bound \eqref{eqn:vorT4b} we use that \eqref{lowbvor} and the trivial relation $r_o(K) \ge h_u(K),\,u\in \mathbb{S}^{d-1},$ imply that $r_o(K)^d\ge d\tau \Sigma(K)^{d/k}$. Since we have assumed \eqref{PVTass1}, this gives for $\Sigma(K)>v$ and $\vt(K)<\eps$ that $\rho_o(K)^d>h(\eps)^dd\tau v^{d/k}$. Hence,  \eqref{eqn:vorT4b} is bounded by
		\begin{align*}
			& \gamma^2 \int_{W} \int_{W} \mathbf 1\{|x-y|\le 2b_\gamma\} \, \mathbb{E}\mathbf 1\{\rho_o(C(y,\eta_\gamma+\delta_{(x,y)}))^d>h(\eps)^d d\tau v_{c,\gamma}^{d/k}\}\nonumber\\
			&\quad \quad \times \mathbf 1\{\Sigma(C(x,\eta_\gamma+\delta_{(x,y)}))>v_{c,\gamma},\,\vt(C(x,\eta_\gamma+\delta_{(x,y)}))<\varepsilon\}\,\mathrm{d}y\,\mathrm{d}x. 
		\end{align*}
		Now we use Lemma \ref{stabaug} and exploit that $\rho_o(C(y,\eta_\gamma+\delta_{(x,y)}))>s$ implies that $\eta_\gamma \cap B_{2s}(y)=\emptyset$. This gives the bound
		\begin{align*}
			&\gamma^2 \int_{W} \int_{W} \mathbf 1\{|x-y|\le 2b_\gamma\}  \, \mathbb{E} \mathbf 1\{\eta_\gamma \cap B_{2h(\eps)(d\tau)^{1/d} v_{c,\gamma}^{1/k}}(y)=\emptyset\}\nonumber\\
			&\quad \quad \times\,\mathbf 1\{\Sigma(C(x,\eta_\gamma\cap K_a(x,y)^c+\delta_{(x,y)}))>v_{c,\gamma},\,\vt(C(y,\eta\cap K_a(x,y)^c+\delta_{(x,y)}))<\varepsilon\}\,\mathrm{d}y\,\mathrm{d}x
		\end{align*}
		where $a:=\frac{h(\eps)}{\sqrt{1-h(\eps)^2}}$. Since the processes $\eta_\gamma \cap K_a(x,y)$ and $\eta_\gamma \cap K_a(x,y)^c$ are independent, we arrive at the bound
		\begin{align*}
			&\gamma^2 \int_{W} \int_{W}  \mathbf 1\{|x-y|\le 2b_\gamma\}  \,\mathbb{P}(\eta_\gamma \cap K_a(x,y) \cap B_{2h(\eps)(d\tau)^{1/d} v_{c,\gamma}^{1/k}}(y)=\emptyset)\nonumber\\
			&\quad \quad \times \,\mathbb{P}(\Sigma(C(x,\eta_\gamma\cap K_a(x,y)^c+\delta_{(x,y)}))>v_{c,\gamma})\,\mathrm{d}y\,\mathrm{d}x.
		\end{align*}
		Now we use Lemma \ref{stabaug} again and note that $\lambda_d( K_a(x,y) \cap B_{r}(y))=r^d \kappa_d k(h(\eps))$ for $r>0$. Hence, the above is bounded by
		\begin{align*}
			&\gamma c\,\lambda_d(W) \lambda_d(B_{2b_\gamma}) \,\exp\Big(-2^{d}  h(\eps)^d d \tau v_{c,\gamma}^{d/k}\kappa_dk(h(\eps))\gamma \Big)\\
			&\quad \le \alpha_6 \log (\gamma) \exp\Big(-2^{d} h(\eps)^d d \tau v_{c,\gamma}^{d/k}\kappa_dk(h(\eps))\gamma \Big). 
		\end{align*}
		
		Finally, we complete the proof and collect the bounds of the $T$-terms. This gives
		\begin{align*}
			&	\mathbf{d_{TV}}(\xi_{c,\gamma} \cap W,\nu_{c,\gamma} \cap W)\le \frac{\alpha_1}{\gamma^2}+\frac{\alpha_2 \log \gamma}{\gamma}+ \frac{\alpha_3}{\gamma}+\alpha_4 \log(\gamma) \exp(-c_0 f(\eps)v_{c,\gamma}^{d/k}\gamma)\\
			&\qquad	+\alpha_5 \gamma^{2-4^d3I\kappa_d(2-\log 2)} \sqrt{\log \gamma}+ \alpha_6 \log (\gamma) \exp\Big(-2^{d} h(\eps)^d d \tau v_{c,\gamma}^{d/k}\kappa_dk(h(\eps))\gamma \Big). 
		\end{align*}
		Using that $\gamma^{-1} v_{c,\gamma}^{-d/k}\log \gamma\to 2^d d\kappa_d \tau$ as $\gamma \to \infty$ we conlude that for all $\delta>0$,
		\begin{align*}
			&	\mathbf{d_{TV}}(\xi_{c,\gamma} \cap W,\nu_{c,\gamma} \cap W)\le C \gamma^{\delta-\min\big[c_0 f(\eps), h(\eps)^d  k(h(\eps))\big]}
		\end{align*}
		for a constant $C>0$ that does not depend on $\gamma$.
	\end{proof}

	\section{Maximum cells in the Poisson-Delaunay mosaic}
	In this section we apply Theorem \ref{maithmallg} to processes of centres of large cells in the Poisson-Delaunay mosaic. As in the previous section we work in the Euclidean space $\mathbb{X}=\mathbb{R}^d$ ($d \ge 2$) with Borel $\sigma$-field $\mathcal{B}^d$ and $d$-dimensional Lebesgue measure $\L_d$. Let $\mu \in \mathbf{N}_s$ and $\mathbf{x}:=(x_1,\dots,x_{d+1})\in \mu^{(d+1)}$ be in general position. Let $B(\mathbf{x})$ be the (unique) open $d$-dimensional ball that has the points $x_1,\dots,x_{d+1}$ on its boundary and let $z(\mathbf{x})$ be its centre. If $\mu \cap B(\mathbf{x})=\emptyset$, we call the simplex $S(\mathbf{x}):=\mathrm{conv}(x_1,\dots,x_{d+1})$ a {\em Delaunay cell}. The system of all such cells is called {\em Delaunay mosaic}.  
	
	Let $\Delta$ be the space of all $d$-simplices with circumcentre at the origin $o$, equipped with the Hausdorff metric. For $k>0$ we call $\Sigma: \Delta \to \mathbb{R}$ a {\em size functional} if it is continuous, $k$-homogeneous and such that $\Sigma$ attains a maximum on the set of simplices with vertices on the unit sphere and if $V_d/\Sigma^{1/k}$ is bounded (where $V_d$ is the volume). Examples for size functionals are the volume and the inradius (see Example \ref{PDTEx1}). We slightly abuse the notation and write $\Sigma(S(\mathbf x)):=\Sigma(S(\mathbf x)-z(\mathbf x))$ for $\mathbf x:=(x_1,\dots,x_{d+1})\in (\mathbb{R}^d)^{d+1}$ in general position. For $c >0$, a threshold $v_{c,\gamma}>0,\,\gamma>0,$ (to be specified below) and a stationary Poisson process $\eta_\gamma$ of intensity $\gamma$ we consider the process
	\begin{align}
		\xi_{c,\gamma}:=\frac{1}{(d+1)!} \sum\limits_{\mathbf{x} \in \eta_\gamma^{(d+1)}} \mathds{1}\{\eta \cap B(\mathbf{x})=\varnothing\}\,\mathds{1}\{\Sigma(S(\mathbf{x}))>v_{c,\gamma}\} \,\delta_{z(\mathbf{x})}.\label{PDTxit}
	\end{align}
	
	Particularly useful in the asymptotic study of random cells is the notion of the typical cell $Z_\gamma$ in a Delaunay mosaic generated by a Poisson process with intensity $\gamma$. This is any random simplex with distribution given by
	\begin{align}
		\P(Z_\gamma \in \cdot)=\frac{1}{\gamma \beta_d (d+1)! } \E \int \mathds{1}\{z(\mathbf{x}) \in [0,1]^d \} \mathds{1}\{S(\mathbf{x})\in \cdot\} \mathds{1}\{\eta_\gamma \cap B(\mathbf{x})=\emptyset\} \eta_\gamma^{(d+1)}(\mathrm{d}\mathbf{x}), \label{PDTtypZ}
	\end{align}
	with some constant $\beta_d>0$ that only depends on $d$ (see \cite{schneider2008stochastic}, p. 450 and (10.31)). This allows us to write the intensity measure of $\xi_{c,\gamma}$ as
	\begin{align}
		\E \xi_{c,\gamma}(A)=\gamma \beta_d \L_d(A)\P(\Sigma(Z_\gamma)>v_{c,\gamma}),\quad A \in \mathcal B^d. \label{PDTintxittyp}
	\end{align}
	We choose the threshold $v_{c,\gamma}$ in \eqref{PDTxit} as
	\begin{align}
		v_{c,\gamma}:=\inf\{v>0:\,\gamma\beta_d\P(\Sigma(Z_\gamma)>v)\le c\}. \label{PDTvtinf}
	\end{align}
	By \eqref{PDTintxittyp} and Lemma \ref{PDTcont} below it holds that $\mathbb{E} \xi_{c,\gamma}=c \lambda_d$.

	For $S \in \Delta$ let $r(S)$ denote the circumradius of $S$ and define $\tau:=\max \{\Sigma(S):\,S \in \Delta,\, r(S)=1\}$. Since $\Sigma$ is $k$-homogeneous we obtain that
	\begin{align} 
		\Sigma(S)\le \tau r(S)^k,\quad S \in \Delta.\label{lowb}
	\end{align}
	If equality holds in \eqref{lowb}, we call $S \in \Delta$ an {\em extremal simplex}. Following \cite{hug2005large} we define a functional $\vt$ that measures the deviation of a simplex $S \in \Delta$ from a regular simplex as follows. Let $S\in \Delta$ and $u_1,\dots,u_{d+1}\in \mathbb{S}^{d-1}$ be such that  $\text{conv}(u_1,\dots,u_{d+1})$ is a regular simplex. Then we define $\vt(S)$ as the smallest number $\alpha>0$ such that there are points $v_1,\dots,v_{d+1}\in \mathbb{S}^{d-1}$ such that $\text{conv}(v_1,\dots,v_{d+1})$ is similar to $S$ and $\|u_i-v_i\| \le \alpha$ for $i\in [d+1]$. Note that $\vt(S)=0$ if and only if $S$ is a regular simplex.
	
	As in \cite{hug2005large} we call $f:[0,1) \to [0,1)$ a {\em stability function} of $\Sigma$ and $\vt$ if $f(0)=0,\,f(\eps) >0$ for $\eps >0$ and 
	\begin{align*}
		\Sigma(S) \le (1-f(\eps) )r(S)^k\tau \quad \text{for all }S \in \Delta\text{ with }\vt(S) \ge \eps.
	\end{align*}
	
	Let $\Sigma$, $\vt$ and $f$ be as above and $v,\,\eps>0$. It was shown in Theorem 1 in \cite{hug2005large} that there is a constant $c_0>0$ depending only on $\Sigma,\,\vt,\,f$ and $d$ such that
	\begin{align}
		\P(\vartheta(Z_\gamma) \ge \eps \mid\,\Sigma(Z_\gamma)\ge v) \le c_1\,\exp(-c_0f(\eps)v^{d/k}\gamma),\quad \gamma>0, \label{af}
	\end{align}
	where $c_1>0$ depends only on $d,\,\eps,\,\Sigma,\,\vt,\,f$.
	
	Next we determine the asymptotic behaviour of $v_{c,\gamma}$ for fixed $c$ as $\gamma \to \infty$. By \eqref{PDTvtinf} and Lemma \ref{PDTcont} we have that $v_{c,\gamma} \to \infty$ as $\gamma \to \infty$. Hence, Theorem 2 in \cite{hug2005large} implies that 
	\begin{align}
		\lim\limits_{\gamma\to \infty} \gamma^{-1} v_{c,\gamma}^{-d/k}\log \P(\Sigma(Z_\gamma)>v_{c,\gamma}) =-\kappa_d \tau^{-d/k}. \label{PDTvass}
	\end{align}
	Since $\gamma^{-1} v_{c,\gamma}^{-d/k} \log \gamma$ is given by
	\begin{align*}
		\gamma^{-1}v_{c,\gamma}^{-d/k}	\log \P(\Sigma(Z_\gamma)>v_{c,\gamma}) \Big(\frac{\log(\gamma\P(\Sigma(Z_\gamma)>v_{c,\gamma}))}{\log \gamma}-1\Big)^{-1}
	\end{align*}
	we conclude from \eqref{PDTvass} and the definition of $v_{c,\gamma}$  that $\gamma^{-1} v_{c,\gamma}^{-d/k} \log \gamma\to \kappa_d \tau^{-d/k}$ as $\gamma \to \infty$.
	
	The next statement is the main result of this section.
	\begin{theorem} \label{PDTmaith}
		Suppose that $\eta_\gamma$ be a stationary Poisson process with intensity $\gamma>1$, $\Sigma$ is a size functional, $\vt$ is a deviation functional with stability function $f$ and $\tau$ from \eqref{lowb}. Assume that all extremal simplices of $\Sigma$ are regular. Let $c>0$, $W \subset \mathbb{R}^d$ be compact and let $\nu_c$ be a stationary Poisson process with intenisty $c$. Then we have for all $\delta>0$ and all $\eps\in (0,1/d)$
		\begin{align*}
			\mathbf{d_{TV}}(\xi_{c,\gamma} \cap W, \nu_c \cap W)\le C \gamma^{\delta-\min[c_0f(\eps),1-I(\eps)]},\quad \gamma>1,
		\end{align*}
		where the constant $C>0$ does not depend on $\gamma$. Here, $I(\eps):=\int_{\frac{d-1}{d}+\eps}^1 (1-s^2)^{d/2}\,\mathrm{d}s$.
	\end{theorem}
	In the following example we discuss applications of Theorem \ref{PDTmaith} to concrete size functionals.
	
	\begin{example} \label{PDTEx1} (a) For the volume $\Sigma=V_d$ it was shown in \cite{hug2004large} that the extremal simplices are regular and that there is a constant $c_d>0$ only depending on $d$ such that for every $\eps\in (0,1)$ we have that
		\begin{align*}
			V_d(S)\le (1-c_d\eps^2)r(S)^d\tau\quad \text{for all }S \in \Delta \text{ with }\vt(S)>\eps.
		\end{align*}
		Hence, a stability function is given by $f(\eps)=c_d \eps^2$. In dimension $d=2$ the distribution of the volume of the typical cell in the Poisson-Delaunay mosaic is known explicitly from \cite{rathie1992volume} and given by
		\begin{align*}
			\P(V_2(Z_\gamma)>v)=\frac{8\pi}{9} \int_{\gamma v}^\infty u\, K_{1/6}^2\left(\frac{2\pi u}{3 \sqrt{3}}\right)\,\mathrm{d}u,\quad v>0,
		\end{align*}
		where $K_{1/6}(u)$ is the modified Bessel function of order $1/6$. Since $K_{1/6}(u)=\sqrt{\frac{\pi}{2u}}e^{u}(1+o(1))$ as $u \to \infty$ we find for $v_{c,\gamma}:=\frac{3 \sqrt{3}}{4 \pi \gamma}\log \big(\frac{3\gamma}{c}\big)$ that
		\begin{align*}
			\lim_{\gamma \to \infty}	\gamma \mathbb P(V_2(Z_\gamma)>v_{c,\gamma})=c.
		\end{align*}
		Letting $c:=e^{-t}$ for $t \in \mathbb{R}$ we obtain from Theorem \ref{PDTmaith}  for all $\eps\in (0,1/2)$:
		\begin{align*}
			\Big|\mathbb{P}\Big(\frac{4 \pi \gamma}{3 \sqrt{3}}\max_{\mathbf x \in \eta_\gamma^{(3)}}\{V_2(S(\mathbf x)):\,z(\mathbf x)\in W, \eta_\gamma \cap B(\mathbf x)=\emptyset\}-\log(3\gamma)\le t\Big)-\exp(-\lambda_2(W)e^{-t})\Big|\le C\gamma^{\delta-M}
		\end{align*}
		for all $t \in \mathbb{R}$ and $\gamma>1$, where $M:=\min[c_0 c_d \eps^2,1/6-\eps+(1/2-\eps)^3/3]$. This quantifies the result from Section 3.2 in \cite{chenavier2014general}.\\
		(b) Let $\Sigma(S)=\rho(S)$ be the inradius (i.e. the radius of the largest ball inscribed in $S$). Then the extremal simplices are also precisely all regular simplices. In was shown in \cite[Section 4]{hug2005large} that for the constant $c_d$ from (a) and $\eps\in [0,1]$ it holds that
		\begin{align*}
			\rho(S)\le (1-c_d\eps^2/d)r(S)\tau\quad \text{for all }S \in \Delta \text{ with }\vt(S)>\eps.
		\end{align*}
		Hence, a stability function is given by $f(\eps)=c_d \eps^2/d$.
	\end{example}
	
	The following spherical Blaschke-Petkantschin formula (Theorem 7.3.1 in \cite{schneider2008stochastic}) will play a fundamental role in the proof of Theorem \ref{PDTmaith}. For $\mathbf{u}:=(u_1,\dots,u_{d+1})$ let $\Delta_{d}(\mathbf{u})$ be the $d$-dimensional volume of the convex hull $\mathrm{conv}(u_1,\dots,u_{d+1})$ of $u_1,\dots,u_{d+1}$. Let $\sigma$ be the uniform probability distribution on the unit sphere $\mathbb{S}^{d-1}$ and $\kappa_d:=\pi^{d/2}/\Gamma(1+d/2),\,d \in \mathbb{N},$ be the volume of the $d$-dimensional unit ball, where $\Gamma(\cdot)$ denotes the Gamma function. Let $f:(\R^d)^{d+1}\to [0,\infty)$ be a measurable function. Then we have
	\begin{align}
		\int\limits_{(\R^d)^{d+1}} f\,\mathrm{d}\lambda_d^{d+1}=d!(d\kappa_d)^{d+1} \int\limits_{\R^d} \int\limits_0^\infty\int\limits_{(\mathbb{S}^{d-1})^{d+1}} f(z+r\mathbf{u}) r^{d^2-1}\,\Delta_d(\mathbf{u})\,\sigma^{d+1}(\mathrm{d}\mathbf{u}) \,\mathrm{d}r\,\mathrm{d}z \label{bpclas}
	\end{align}
	where $z+r \mathbf{u}:=(z+ru_1,\dots,z+ru_{d+1})$.
	
	From the Mecke formula and \eqref{PDTtypZ} we obtain for all $v \ge 0$ that
	\begin{align*}
		\mathbb{P}(r(Z_\gamma)>v)=\frac{\gamma^{d+1}}{\beta_d(d+1)!} \int \mathbf 1\{z(\mathbf x) \in [0,1]^d\} \mathbf 1\{r(\mathbf x)>v\} e^{-\gamma \kappa_d r(\mathbf x)^d}\,\lambda_d^{d+1}(\mathrm{d}\mathbf x),
	\end{align*}
	which is by \eqref{bpclas} given by
	\begin{align}
		\frac{\gamma \kappa_d d^dC_d}{\beta_d(d+1)} \int_{(\frac{v}{\gamma \kappa_d})^{1/d}}^\infty e^{-s} s^{d-1} \mathrm d s\quad \text{with}\quad C_d:=\int_{(\mathbb S^{d-1})^{d+1}} \Delta_d(\mathbf u)\,\sigma_{d-1}^{d+1}(\mathrm{d} \mathbf u).\label{PDT:defCd}
	\end{align}
	Using here that $\mathbb{P}(r(Z_\gamma)>0)=1$ and that $\int_0^\infty e^{-s} s^{d-1}\,\mathrm d s=\Gamma(d)$, we find
	\begin{align}
		\mathbb{P}(\gamma \kappa_d r(Z_\gamma)^d>v)=\frac{1}{\Gamma(d)} \int_v^\infty e^{-s} s^{d-1}\,\mathrm{d}s,\quad v \ge 0,\label{PDTdistrd}
	\end{align}
	i.e. $\gamma \kappa_d r(Z_\gamma)^d$ follows a Gamma distribution with parameter $d$.
	
	\begin{figure}[bt]
		\centering
		\includegraphics[scale=0.79]{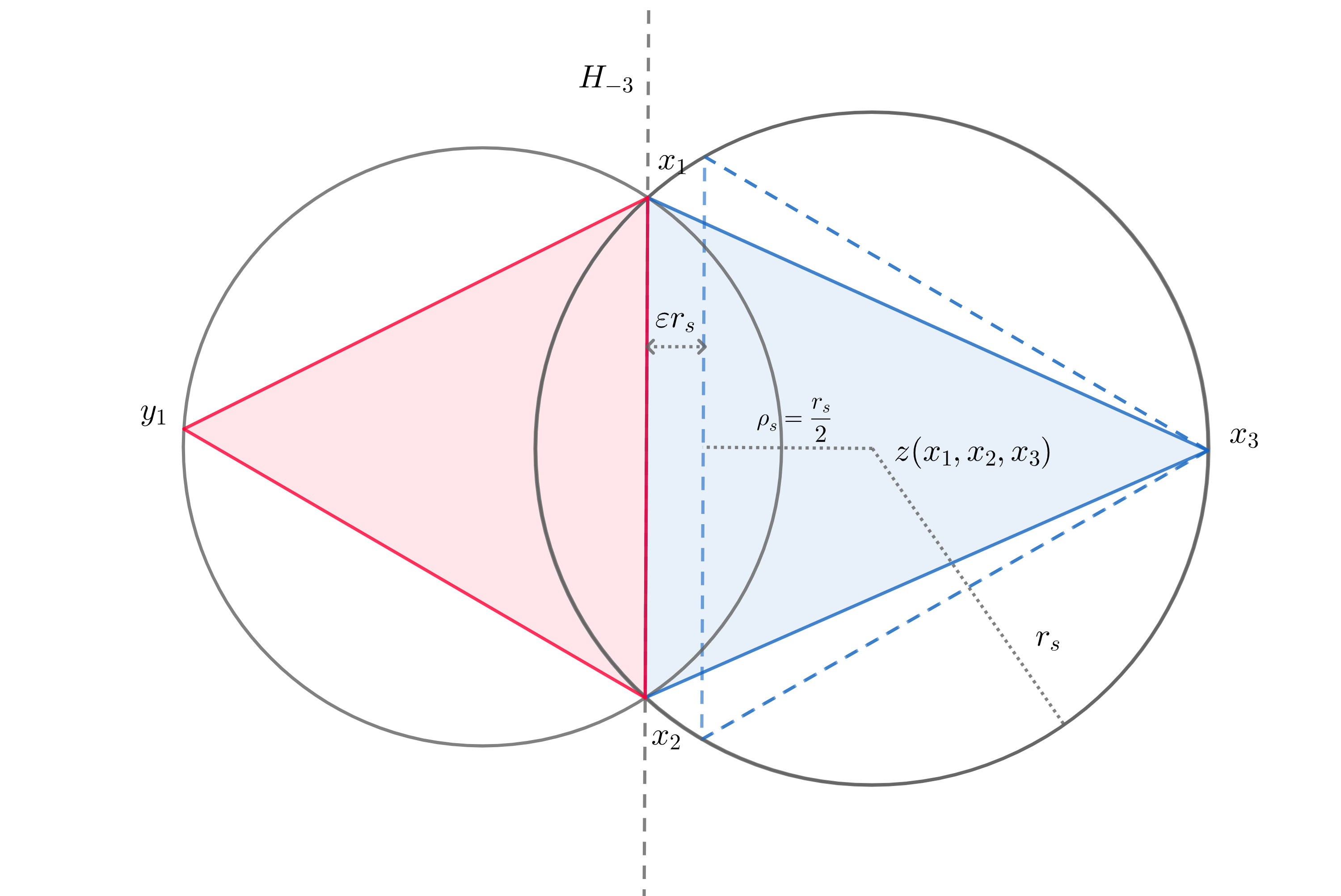}
		\caption{The figure shows the situation of Lemma \ref{PDTform} for $d=2$ and $\ell=2$. The dashed drawn triangle is regular.}
	\end{figure}
	\begin{lemma} \label{PDTform}
		For $\ell \in [d+1]$ let  $\mathbf{x}\in (\R^d)^{d+1}$ and $\mathbf{y}\in (\R^d)^{\ell}$ such that $\mathbf{x}$ and $(\mathbf{x}_\ell,\mathbf{y})$ are in general position. We assume that 
		\begin{align*}
			B(\mathbf{x})\cap \{y_1,\dots,y_{d+1-\ell}\}=\emptyset,\qquad
			B(\mathbf{x}_\ell,\mathbf{y}) \cap \{x_{\ell+1},\dots,x_{d+1}\}=\emptyset.
		\end{align*}
		and that $\vt(\mathbf x)<\eps$ and $\vt(\mathbf{x}_\ell,\mathbf y)<\eps$ for some $\eps \in (0,1/d)$.
		Then we have
		\begin{align*}
			\lambda_d(B(\mathbf{x})\cup B(\mathbf{x}_\ell,\mathbf{y}))\ge \big(1-I(\eps)\big) \kappa_d r(\mathbf{x})^d+\kappa_d r(\mathbf{x}_\ell,\mathbf{y})^d
		\end{align*}
		with $I(\eps):=\int_{\frac{d-1}{d}+\eps}^1 (1-s^2)^{d/2}\,\mathrm{d}s$.
	\end{lemma}
	
	\begin{proof} Let $\rho_s$ be the inradius (radius of the largest inscribed sphere) and $r_s$ be the circumradius (radius of the smallest circumscribing sphere) of a regular $(d+1)$-simplex. We use the fact that $d\rho_s=r_s$. For $1 \le k \le d+1$ let $H_{-k}$ be the (unique) hyperplane through $\{x_1,\dots,x_{d+1}\} \setminus \{x_k\}$. It follows from the definition of $\vt$ that $\vt(\mathbf{x})<\eps$ implies that 
		\begin{align}
			\min_{1 \le k \le d+1} d(z(\mathbf{x}),H_{-k}) \le r(\mathbf{x})\Big(\frac{1}{d}+\eps\Big). \label{eqn:mindist}
		\end{align}
		Let $H_{-k}^-$ be the closed halfspace that is bounded by $H_{-k}$ and that does not contain $z(\mathbf x)$. We have
		\begin{align}
			\lambda_d(B(\mathbf{x}) \cup B(\mathbf{x}_\ell,\mathbf{y}))\ge \lambda_d(B(\mathbf{x}))+\lambda_d(B(\mathbf{x}_\ell,\mathbf{y}))-\max_{1 \le k \le d+1} \lambda_d(B(\mathbf{x}) \cap H_{-k}^-). \label{eqn:volcup}
		\end{align}
		For  $d(z(\mathbf{x}),H_{-k})=ar(\mathbf{x})$ we have 
		\begin{align*}
			\lambda_d(B(\mathbf{x}) \cap H_{-k}^-)=r(\mathbf{x})^d\kappa_d\int_{1-a}^1 (1-s^2)^{d/2}\,\mathrm{d}s.
		\end{align*}
		Hence, we find the assertion from \eqref{eqn:mindist} and \eqref{eqn:volcup}.
	\end{proof}

	\begin{proof}[Proof of Theorem \ref{PDTmaith}]
		Let $m=d+1$, $\gamma>0$ and $c>0$ and $W \subset \mathbb{R}^d$ be compact. We apply Theorem \ref{maithmallg} with $\lambda:=\gamma \lambda_d\cap W$ and
		\begin{align*}
			g(\mathbf{x},\mu)=\mathbf{1}\{z(\mathbf{x}) \in W\}\mathbf{1}\{\mu \cap B(\mathbf{x})=\varnothing,\,\Sigma(S(\mathbf{x}))>v_{c,\gamma}\}.
		\end{align*}
		We choose the stabilization radius $R(z,\mu):=\|z-x_1\|$ if there is a unique $(d+1)$-tuple $\mathbf{x}=(x_1,\dots,x_{d+1})\in \mu^{(d+1)}$ (up to permutations of the components of $\mathbf{x}$) such that $z(\mathbf{x})=z$ (in this case $R(z,\mu)=r(\mathbf x)$) and $R(z,\mu):=\infty$, otherwise. Let
		\begin{align}
			b_{\mathbf x}:= b_\gamma=(3\gamma^{-1} \kappa_d^{-1}\log \gamma)^{1/d}. \label{PDTdefbt}
		\end{align}
		This gives the bound
		\begin{align*}
			\mathbf{d_{TV}}(\xi_{c,\gamma} \cap W, \nu_c \cap W)\le \|\mathbb{E}(\xi_{c,\gamma}\cap W)-\mathbb{E}(\nu_{c}\cap W)\|+T_1+T_2+T_3+T_4+T_5
		\end{align*}
		where
		\begin{align*}
			T_1&=\frac{2\gamma^{d+1}\mathbb{E}\xi_{n,c}(W)}{(d+1)!} \int\mathbf 1 \{z(\mathbf x)\in W\} \mathds 1\{\Sigma(S(\mathbf x))>v_{c,\gamma}\}  \mathbf 1\{r(\mathbf x)>b_\gamma\}\mathbb{P}(\eta_\gamma \cap B(\mathbf x)=\emptyset)\,\lambda_d^{d+1}(\mathrm{d}\mathbf x),\\
			T_2&=\frac{\gamma^{2d+2}}{((d+1)!)^2} \iint \mathbf 1 \{z(\mathbf x)\in W, |z(\mathbf{x})-z(\mathbf y)| \le 2b_\gamma\} \mathbf 1\{\Sigma(S(\mathbf x))>v_{c,\gamma},\Sigma(S(\mathbf y))>v_{c,\gamma}\} \\
			&\qquad \qquad \times  \mathbb{P}(\eta_\gamma \cap B(\mathbf x)=\emptyset)\, \mathbb{P}(\eta_\gamma \cap B(\mathbf y)=\emptyset)\,\lambda_d^{d+1}(\mathrm{d}\mathbf y)\,\lambda_d^{d+1}(\mathrm{d}\mathbf x),\\
			T_3&=\frac{\gamma^{2d+2}}{((d+1)!)^2} \iint \mathbf 1 \{z(\mathbf x)\in W, z(\mathbf{y})\in W\}\mathbf 1\{\Sigma(S(\mathbf x))>v_{c,\gamma},\Sigma(S(\mathbf y))>v_{c,\gamma}\} \\
			&\qquad \qquad \times  \mathbf 1\{r(\mathbf x)>b_\gamma\}\mathbb{P}((\eta_\gamma+\delta_{(\mathbf x,\mathbf y)}) \cap (B(\mathbf x) \cup B(\mathbf y))=\emptyset)\,\lambda_d^{d+1}(\mathrm{d}\mathbf y)\,\lambda_d^{d+1}(\mathrm{d}\mathbf x),\\
			T_4&=\frac{2\gamma^{2d+2}}{((d+1)!)^2} \iint \mathbf 1 \{z(\mathbf x)\in W, |z(\mathbf{x})-z(\mathbf y)| \le 2b_\gamma\}\mathbf 1\{\Sigma(S(\mathbf x))>v_{c,\gamma},\Sigma(S(\mathbf y))>v_{c,\gamma}\} \\
			&\qquad \qquad \times  \mathbf 1\{r(\mathbf x)>b_\gamma\}\mathbb{P}((\eta_\gamma+\delta_{(\mathbf x,\mathbf y)}) \cap (B(\mathbf x) \cup B(\mathbf y))=\emptyset)\,\lambda_d^{d+1}(\mathrm{d}\mathbf y)\,\lambda_d^{d+1}(\mathrm{d}\mathbf x),\\
			T_5&=\sum_{\ell=1}^{d}\frac{\gamma^{2d+2-\ell}}{\ell!(d+1-\ell)!} \iint \mathbf 1 \{z(\mathbf x)\in W, z(\mathbf{x}_\ell,\mathbf y)\in W\}\mathbf 1\{\Sigma(S(\mathbf x))>v_{c,\gamma},\Sigma(S(\mathbf x_\ell,\mathbf y))>v_{c,\gamma}\} \\
			&\qquad \qquad \times \mathbb{P}((\eta_\gamma+\delta_{(\mathbf x,\mathbf y)}) \cap (B(\mathbf x) \cup B(\mathbf x_\ell,\mathbf y))=\emptyset)\,\lambda_d^{d+1-\ell}(\mathrm{d}\mathbf y)\,\lambda_d^{d+1}(\mathrm{d}\mathbf x).
		\end{align*}
		From \eqref{PDTintxittyp} and the choice of $v_{c,\gamma}$ in \eqref{PDTvtinf} we find that $\mathbb{E} \xi_{c,\gamma}=\mathbb{E} \nu_c$. Hence, $\|\mathbb{E}(\xi_{c,\gamma}\cap W)-\mathbb{E}(\nu_{c}\cap W)\|=0$, $\gamma>1$. Next we bound the $T$-terms from the right-hand side. In the following, $\alpha_i>0$ ($i \in \mathbb{N}$) are positive constants that do not depend on $\gamma$. Their precise values are not important for the argument.
		
		{\em The estimate of $T_1$.} From the definition of $b_\gamma$ we have that $T_1$ is bounded by
		\begin{align*}
			\frac{2\gamma^{d+1} c \lambda_d(W)}{(d+1)!} \int\mathbf 1 \{z(\mathbf x)\in W\}   \mathbf 1\{\gamma \kappa_dr(\mathbf x)^d>3 \log \gamma\} e^{-\gamma \kappa_d r(\mathbf{x})^d}\,\lambda_d^{d+1}(\mathrm{d}\mathbf x).
		\end{align*}
		Now we invoke the definition of the typical cell $Z_\gamma$ from \eqref{PDTtypZ} and use \eqref{PDTdistrd}. This gives for the above
		\begin{align}
			2c\gamma \lambda_d(W)^2\beta_d\mathbb{P}(\gamma \kappa_d r(Z_\gamma)^d>3 \log \gamma)=\frac{2c\gamma \lambda_d(W)^2\beta_d }{\Gamma(d)} \int_{3 \log \gamma}^\infty e^{-s} s^{d-1} \,\mathrm ds\le \frac{\alpha_1 (\log \gamma)^{d-1}}{\gamma^2}.\label{PDT1bou}
		\end{align}
		
		{\em The estimate of $T_2$.} By definition of the typical cell $Z_\gamma$  we find that $T_2$ is given by
		\begin{align}
			\gamma^2 \lambda_d(B_{2b_\gamma}) \lambda_d(W) \beta_d^2 \mathbb{P}(\Sigma(Z_\gamma)>v_{c,\gamma})^2 \le \frac{\alpha_2 \log \gamma}{\gamma}. \label{PDT2bou} 
		\end{align}
		
		{\em The estimate of $T_3$.} For the estimate of $T_3$ we assume (at the cost of a factor 2) that $r(\mathbf y)\le r( \mathbf x)$ and bound $T_3$ by
		\begin{align*}
			&\frac{\gamma^{2d+2}}{((d+1)!)^2} \iint \mathbf 1 \{z(\mathbf x)\in W, z(\mathbf{y})\in W\}  \mathbf 1\{r(\mathbf x)>b_\gamma\} \mathbf 1\{r(\mathbf x)\ge r(\mathbf y)\} \\
			&\qquad \qquad \times \mathbb{P}(\eta_\gamma\cap (B(\mathbf x) \cup B(\mathbf y)) =\emptyset)\,\lambda_d^{d+1}(\mathrm{d}\mathbf y)\,\lambda_d^{d+1}(\mathrm{d}\mathbf x).
		\end{align*}
		Now we apply \eqref{bpclas} to the inner integral and obtain for the above
		\begin{align*}
			&\frac{\gamma^{2d+2} \kappa_d d^d\lambda_d(W)C_d d!}{d^2((d+1)!)^2} \int r(\mathbf x)^{d^2} \mathbf 1 \{z(\mathbf x)\in W\}\mathbf 1\{r(\mathbf x)>b_\gamma\}  \mathbb{P}(\eta_\gamma\cap B(\mathbf x) =\emptyset)\,\lambda_d^{d+1}(\mathrm{d}\mathbf x)
		\end{align*}
		with the constant $C_d$ from \eqref{PDT:defCd}. We apply \eqref{bpclas} a second time and substitute $s:=\gamma \kappa_d r^d$. Thus we arrive at
		\begin{align*}
			\alpha_3 \gamma^{2d+2}  \int_{b_\gamma}^\infty r^{2d^2-1} e^{-\gamma \kappa_d r^d} \mathrm d r=\alpha_4 \gamma^2 \int_{\gamma \kappa_d b_\gamma^d}^\infty s^{2d-1} e^{-s} \,\mathrm ds\le \frac{\alpha_5 (\log \gamma)^{2d-1}}{\gamma} .
		\end{align*}
		
		{\em The estimate of $T_4$.} Let $\eps>0$. We consider the shapes of the simplices $S(\mathbf x)$ and $S(\mathbf y)$ and split $T_4$ into
		\begin{align}
			&\frac{2\gamma^{2d+2}}{((d+1)!)^2} \iint \mathbf 1\{\vt(S(\mathbf x))\le \eps ,\,\vt(S(\mathbf y))\le \eps \}\mathbf 1 \{z(\mathbf x)\in W, |z(\mathbf{x})-z(\mathbf y)| \le 2b_\gamma,\,r(\mathbf x)\le b_\gamma,\,r(\mathbf y)\le b_\gamma\}\nonumber\\
			&\qquad \times \mathbf 1\{\Sigma(S(\mathbf x))>v_{c,\gamma},\Sigma(S(\mathbf y))>v_{c,\gamma}\} \mathbb{P}((\eta_\gamma+\delta_{(\mathbf x,\mathbf y)}) \cap (B(\mathbf x) \cup B(\mathbf y))=\emptyset)\,\lambda_d^{d+1}(\mathrm{d}\mathbf y)\,\lambda_d^{d+1}(\mathrm{d}\mathbf x)\label{PDT4bouleps}\\
			&\quad  + \frac{4\gamma^{2d+2}}{((d+1)!)^2} \iint \mathbf 1\{\vt(S(\mathbf x))>\eps\}\mathbf 1 \{z(\mathbf x)\in W, |z(\mathbf{x})-z(\mathbf y)| \le 2b_\gamma\}\mathbf 1\{r(\mathbf x)\le b_\gamma,\,r(\mathbf y)\le b_\gamma\} \nonumber\\
			&\qquad\times \mathbf 1\{\Sigma(S(\mathbf x))>v_{c,\gamma},\Sigma(S(\mathbf y))>v_{c,\gamma}\}\mathbb{P}((\eta_\gamma+\delta_{(\mathbf x,\mathbf y)}) \cap (B(\mathbf x) \cup B(\mathbf y))=\emptyset)\,\lambda_d^{d+1}(\mathrm{d}\mathbf y)\,\lambda_d^{d+1}(\mathrm{d}\mathbf x).\label{PDT4bougeps}
		\end{align}
		To estimate \eqref{PDT4bouleps} we assume (at the cost of a factor 2) that $r(\mathbf x)<r(\mathbf y)$ and use that Lemma \ref{PDTform} implies that for $\vt(S(\mathbf x))<\eps$ and $\vt(S(\mathbf y))<\eps$,
		$$
		\lambda_d(B(\mathbf x)\cup B(\mathbf y)) \ge (1-I(\eps))\kappa_d r(\mathbf{x})^d+\kappa_d r(\mathbf y)^d\ge (2-I(\eps)) \kappa_d r(\mathbf {x})^2.
		$$
		where $I(\eps):=\int_{\frac {d-1}{d}+\eps}^1 (1-s^2)^{d/2}\,\mathrm{d}s$. Since $|z(\mathbf x)-z(\mathbf y)|\le 2b_\gamma$ and $r(\mathbf y)\le b_\gamma$ implies by the triangle inequality that $|z(\mathbf x)-y_i| \le 3b_\gamma$, $i=1,\dots,d+1$, we find for \eqref{PDT4bouleps} the bound
		\begin{align}
			\frac{2 \gamma^{2d+2}\lambda_d(B_{3b_\gamma})^{d+1}}{((d+1)!)^2} \int \mathbf 1\{z(\mathbf x)\in W\} \mathbf 1\{\Sigma(S(\mathbf x))>v_{c,\gamma}\} e^{-\gamma \kappa_d (2-I(\eps)) r(\mathbf{x})^d} \,\lambda_d^{d+1}(\mathrm{d}\mathbf{x}).\label{PDT4bouleps1}
		\end{align}
		Now we invoke \eqref{lowb} which says that $\tau r(\mathbf x)^k \ge\Sigma(S(\mathbf x))$ for $\mathbf x$ in general position. Hence, $r(\mathbf x)>(v_{c,\gamma}/\tau)^{1/k}$ for $\Sigma(S(\mathbf x))>v_{c,\gamma}$. Therefore, the above is bounded by
		\begin{align}
			\frac{2 \gamma^{2d+2}\lambda_d(B_{3b_\gamma})^{d+1}e^{-\gamma \kappa_d(1-I(\eps))\tau^{-d/k}v_{c,\gamma}^{d/k}}}{((d+1)!)^2} \int \mathbf 1\{z(\mathbf x)\in W\} \mathbf 1\{\Sigma(S(\mathbf x))>v_{c,\gamma}\} e^{-\gamma \kappa_d r(\mathbf{x})^d} \,\lambda_d^{d+1}(\mathrm{d}\mathbf{x}). \label{PDTvolbou}
		\end{align}
		Since the integral is by definition of $v_{c,\gamma}$ equal to $c \lambda_d(W)$, we find from the definiton of $b_\gamma$ for \eqref{PDTvolbou} the bound 
		\begin{align*}
			\alpha_6 (\log \gamma)^{d+1} \exp(-\gamma \kappa_d(1-I(\eps))\tau^{-d/k}v_{c,\gamma}^{d/k}).
		\end{align*}
		
		An analogous application of the triangle inequality as above yields for \eqref{PDT4bougeps} the bound
		\begin{align*}
			\frac{4 \gamma^{2d+2}\lambda_d(B_{3b_\gamma})^{d+1}}{((d+1)!)^2} \int \mathbf1 \{\vt(S(\mathbf x))>\eps\}\mathbf 1\{z(\mathbf x)\in W\} \mathbf 1\{\Sigma(S(\mathbf x))>v_{c,\gamma}\} e^{-\gamma \kappa_d r(\mathbf{x})^d} \,\lambda_d^{d+1}(\mathrm{d}\mathbf{x}),
		\end{align*}
		for which we obtain by the definition of the typical cell $Z_\gamma$ and by \eqref{af} the bound
		\begin{align*}
			&\alpha_7 \gamma^{2d+2} \lambda_d(B_{3b_n})^{d+1} \mathbb{P}(\Sigma(Z_\gamma)>v_{c,\gamma}) \mathbb{P}(\vt (Z_\gamma)>\eps \mid \, \Sigma(Z_\gamma)>v_{c,\gamma})\\
			&\quad\le \alpha_8 (\log \gamma)^{d+1} \exp(-c_0 f(\eps) v_{c,\gamma}^{d/k}\gamma).
		\end{align*}
		
		{\em The estimate of $T_5$.} To estimate $T_5$ we distinguish by the circumradii $r(\mathbf x)$ and $r(\mathbf x_\ell,\mathbf y)$. This gives
		\begin{align}
			&\sum_{\ell=1}^{d}\frac{\gamma^{2d+2-\ell}}{\ell!(d+1-\ell)!} \iint \mathbf 1\{r(\mathbf x)\le b_\gamma,\,r(\mathbf x_\ell,\mathbf y)\le b_\gamma\}\,\mathbf 1\{\Sigma(S(\mathbf x))>v_{c,\gamma},\Sigma(S(\mathbf x_\ell,\mathbf y))>v_{c,\gamma}\}\ \nonumber \\
			&\quad  \times\mathbf 1 \{z(\mathbf x)\in W, z(\mathbf{x}_\ell,\mathbf y)\in W\} \mathbb{P}((\eta_\gamma+\delta_{(\mathbf x,\mathbf y)}) \cap (B(\mathbf x) \cup B(\mathbf x_\ell,\mathbf y))=\emptyset) \,\lambda_d^{d+1-\ell}(\mathrm{d}\mathbf y)\,\lambda_d^{d+1}(\mathrm{d}\mathbf x)\label{PDT5bourleb}\\
			&\, +\sum_{\ell=1}^{d}\frac{2\gamma^{2d+2-\ell}}{\ell!(d+1-\ell)!} \iint \mathbf 1\{r(\mathbf x)>b_n,\,r(\mathbf x)\ge r(\mathbf x_\ell,\mathbf y)\}\,\mathbf 1\{\Sigma(S(\mathbf x))>v_{c,\gamma},\Sigma(S(\mathbf x_\ell,\mathbf y))>v_{c,\gamma}\} \nonumber\\
			&\quad \times\,\mathbf 1 \{z(\mathbf x)\in W, z(\mathbf{x}_\ell,\mathbf y)\in W\}\mathbb{P}((\eta_\gamma+\delta_{(\mathbf x,\mathbf y)}) \cap (B(\mathbf x) \cup B(\mathbf x_\ell,\mathbf y))=\emptyset)\, \lambda_d^{d+1-\ell}(\mathrm{d}\mathbf y)\,\lambda_d^{d+1}(\mathrm{d}\mathbf x)\label{PDT5bourgeb}
		\end{align}
		where the factor 2 in the second term comes from the assumption $r(\mathbf x)\ge r(\mathbf x_\ell,\mathbf y)$. In \eqref{PDT5bourleb} we next distinguish by the shape of the simplices $S(\mathbf x)$ and $S(\mathbf x_\ell,\mathbf y)$. This gives the bound
		\begin{align}
			&\sum_{\ell=1}^{d}\frac{\gamma^{2d+2-\ell}}{\ell!(d+1-\ell)!} \iint \mathbf 1\{\vt(S(\mathbf x))\le \eps,\,\vt(S(\mathbf x_\ell,\mathbf y))\le \eps\} \mathbf 1\{r(\mathbf x)\le b_\gamma,\,r(\mathbf x_\ell,\mathbf y)\le b_\gamma\}  \nonumber \\
			&\quad \times \mathbf 1 \{z(\mathbf x)\in W, z(\mathbf{x}_\ell,\mathbf y)\in W\}\, \mathbb{P}((\eta_\gamma+\delta_{(\mathbf x,\mathbf y)}) \cap (B(\mathbf x) \cup B(\mathbf x_\ell,\mathbf y))=\emptyset)\,\lambda_d^{d+1-\ell}(\mathrm{d}\mathbf y)\,\lambda_d^{d+1}(\mathrm{d}\mathbf x) \label{PDT5bourlebleps}\\
			&\, +\sum_{\ell=1}^{d}\frac{2\gamma^{2d+2-\ell}}{\ell!(d+1-\ell)!} \iint \mathbf 1\{\vt(S(\mathbf x))>\eps\} \mathbf 1\{r(\mathbf x)\le b_n,\,r(\mathbf x_\ell,\mathbf y)\le b_n\} \mathbf 1\{\Sigma(S(\mathbf x))>v_{c,\gamma}\} \nonumber \\
			&\quad \times\mathbf 1 \{z(\mathbf x)\in W, z(\mathbf{x}_\ell,\mathbf y)\in W\}\, \mathbb{P}((\eta_\gamma+\delta_{(\mathbf x,\mathbf y)}) \cap (B(\mathbf x) \cup B(\mathbf x_\ell,\mathbf y))=\emptyset)\,\lambda_d^{d+1-\ell}(\mathrm{d}\mathbf y)\,\lambda_d^{d+1}(\mathrm{d}\mathbf x) \label{PDT5bourlebgeps}.
		\end{align}
		For \eqref{PDT5bourlebleps} we exploit now the bound \eqref{PDTvolbou} and use that by the triangle inequality $|z(\mathbf x)-y_i|\le |z(\mathbf x)-x_1|+|x_1-z(\mathbf x_\ell,\mathbf y)|+|z(\mathbf x_\ell,\mathbf y)-y_i| \le r(\mathbf x)+2r(\mathbf x_\ell,\mathbf y)$ for $i=1,\dots,d+1-\ell$. This gives the bound
		\begin{align*}
			&\sum_{\ell=1}^{d}\frac{\gamma^{2d+2-\ell}\lambda_d(B_{3b_\gamma})^{d+1-\ell}}{\ell!(d+1-\ell)!} \int \mathbf 1 \{z(\mathbf x)\in W\}\mathbf 1\{\Sigma(S(\mathbf x))>v_{c,\gamma}\} e^{-\gamma \kappa_d (2-I(\eps))r(\mathbf x)^d}\,\lambda_d^{d+1}(\mathrm{d}\mathbf x),
		\end{align*}
		which can be bounded exactly as \eqref{PDT4bouleps1} above and we arrive at the bound
		\begin{align*}
			\alpha_9 (\log \gamma)^{d} \exp(-\gamma \kappa_d(1-I(\eps))\tau^{-d/k}v_{c,\gamma}^{d/k}).
		\end{align*}
		
		For \eqref{PDT5bourlebgeps} we use again that $|z(\mathbf x)-y_i|\le r(\mathbf x)+2r(\mathbf x_\ell,\mathbf y)$ for $i=1,\dots,d+1-\ell$. Analogously to the estimate of \eqref{PDT4bougeps} we find the bound
		\begin{align*}
			\alpha_{10} (\log \gamma)^{d} \exp(-c_0 f(\eps) v_n^{d/k}\gamma).
		\end{align*}
		
		Finally, we dicuss \eqref{PDT5bourgeb}. Since $|z(\mathbf x)-y_i|\le 3r(\mathbf x)$ for $r(\mathbf x_\ell,\mathbf y)\le r(\mathbf x)$, we find the bound
		\begin{align*}
			&\sum_{\ell=1}^{d}\frac{2\gamma^{2d+2-\ell}}{\ell!(d+1-\ell)!} \int \lambda_d(B_{3r(\mathbf x)})^{d+1-\ell}\mathbf 1\{r(\mathbf x)>b_\gamma\} \mathbf 1 \{z(\mathbf x)\in W\} \mathbb{P}(\eta_\gamma \cap B(\mathbf x) =\emptyset)\,\lambda_d^{d+1}(\mathrm{d}\mathbf x)\\
			&\, =\sum_{\ell=1}^{d}\frac{2(3^d\kappa_d)^{d+1-\ell}\gamma^{2d+2-\ell}}{\ell!(d+1-\ell)!} \int r(\mathbf x)^{d(d+1-\ell)}\mathbf 1\{r(\mathbf x)>b_\gamma\} \mathbf 1 \{z(\mathbf x)\in W\} \mathbb{P}(\eta_\gamma \cap B(\mathbf x) =\emptyset)\,\lambda_d^{d+1}(\mathrm{d}\mathbf x).
		\end{align*}
		which is by \eqref{bpclas} bounded by
		\begin{align*}
			d! \omega_d^{d+1} \lambda_d(W) \sum_{\ell=1}^{d}\frac{2(2^d\kappa_d)^{d+1-\ell}\gamma^{2d+2-\ell}}{\ell!(d+1-\ell)!} \int_{b_\gamma}^\infty r^{d(d+1-\ell)+d^2-1 } e^{-\gamma \kappa_dr^d} \,\mathrm{d}r\le \frac { \alpha_{11} (\log \gamma)^{d}} {\gamma}.
		\end{align*}
		
		Finally, we collect the bounds of the different $T$-terms and obtain that 
		\begin{align*}
			\mathbf{d_{TV}}(\xi_{c,\gamma} \cap W, \nu_c \cap W)&\le \frac{\alpha_1 (\log \gamma)^{d-1}}{\gamma^2}+\frac{\alpha_2 \log \gamma}{\gamma}+ \frac{\alpha_5 (\log \gamma)^{2d-1}}{\gamma}+\frac { \alpha_{11} (\log \gamma)^{d}} {\gamma}\\
			&\quad+(\alpha_6 (\log \gamma)^{d+1}+\alpha_9 (\log \gamma)^{d}) \exp(-\gamma \kappa_d(1-I(\eps))\tau^{-d/k}v_{c,\gamma}^{d/k})\\
			&\quad+(\alpha_8 (\log \gamma)^{d+1}+\alpha_{10} (\log \gamma)^{d}) \exp(-c_0 f(\eps) v_n^{d/k}\gamma).
		\end{align*}
		Hence, we conclude from the asyomptotic form of $v_{c,\gamma}$ that for all $\delta>0$ and some constant $C>0$
		\begin{align*}
			\mathbf{d_{TV}}(\xi_{c,\gamma} \cap W, \nu_c \cap W)\le C \gamma^{\delta-\min[c_0f(\eps),1-I(\eps)]}.
		\end{align*}
	\end{proof}
	
	\begin{lemma} \label{PDTcont}
		For all $\gamma>0$ the distribution $\mathbb{P}^{\Sigma(Z_\gamma)}$ of $\Sigma(Z_\gamma)$ and the Lebesgue measure $\lambda_1$ are equivalent on $[0,\infty).$ 
	\end{lemma}
	
	\begin{proof}
		First we show that $\mathbb{P}^{\Sigma(Z_\gamma)}$ is absolutely continuous with respect to $\lambda_1$ on $[\gamma^{-1/k},\infty)$. For $A \in \mathcal B^1$ we have by Theorem 7.3.1 in \cite{schneider2008stochastic} 
		\begin{align}
			\mathbb P(\Sigma(Z_\gamma)\in A)=\frac{\gamma^d}{\beta_d} \int_{(\mathbb S^{d-1})^{d+1}} \int_0^\infty \mathbf 1\{\Sigma(\mathrm{conv}(r\mathbf u))\in A\} e^{-\gamma \kappa_d r^d} r^{d^2-1} \Delta_d(\mathbf u) \,\mathrm d r\,\sigma^{d+1}(\mathrm d \mathbf u). \label{PDTabscont}
		\end{align}
		Since $\Sigma$ is $k$-homogeneous and $\Sigma(\mathrm{conv}(\mathbf u))$ is assumed to be bounded, we have that $\Sigma(\mathrm{conv}(r\mathbf u))\in A$ if $r^k \in A/\Sigma(\mathrm{conv}(\mathbf u))$. Thus, the inner integral in \eqref{PDTabscont} vanishes if $A$ is a $\lambda_1$-null set.
		
		Next we show that the Radon-Nikod\'{y}m density of $\mathbb{P}^{\Sigma(Z_\gamma)}$ with respect to the Lebesgue measure $\lambda_1$ is positive on $[0, \infty)$. From Lemma 1 in \cite{hug2005large} we obtain for $a >0$ and all $\eps>0$
		\begin{align*}
			\frac{\mathrm{d}\mathbb{P}^{\Sigma(Z_\gamma)}}{\mathrm{d}\lambda_1} (a) =\lim_{h \downarrow 0} \frac{\mathbb{P}(\Sigma(Z_\gamma) \in a[1,1+h))}{ah}\ge  c_1(a^{d/k}\gamma)^d \exp\Big(-\frac{\kappa_d}{\tau^{d/k}} (1+\eps) a^{d/k}\gamma\Big)>0.
		\end{align*}
		This gives that $\mathbb{P}^{\Sigma(Z_\gamma)}$ and $\lambda_1$ are equivalent measures on $[0, \infty)$.
	\end{proof}
	
\noindent
{\bf Acknowledgments:} 
		The author wishes to thank Günter Last for helpful discussions.

	
	

\end{document}